\documentclass[12pt]{amsart}
\usepackage{amssymb}
\usepackage{amsmath,amssymb,amsfonts,amsthm,latexsym}
\usepackage{mathrsfs}
\usepackage{color}
\usepackage{eucal}
\usepackage{eufrak}
\usepackage[all]{xypic}
\usepackage{xspace}
\usepackage{hyperref}
\usepackage{natbib}
\usepackage{graphicx}

\theoremstyle{plain}
\newtheorem{theorem}{Theorem}
\newtheorem{corollary}{Corollary}

\newtheorem{lemma}{Lemma}
\newtheorem{proposition}{Proposition}

\theoremstyle{definition}
\newtheorem{definition}{Definition}

\theoremstyle{example}

\theoremstyle{remark}

\numberwithin{equation}{section}

\begin{document}
\begin{center}
{\bf\Large Combinatorics of RNA-RNA interaction}
\\
\vspace{15pt} Thomas J. X. Li, Christian M. Reidys$^{\,\star}$
\end{center}

\begin{center}
Center for Combinatorics, LPMC-TJKLC\\
Nankai University  \\
Tianjin 300071\\
         P.R.~China\\
         Phone: *86-22-2350-6800\\
         Fax:   *86-22-2350-9272\\
duck@santafe.edu
\end{center}

\centerline{\bf Abstract}{\small
RNA-RNA binding is an important phenomenon observed for many classes
of non-coding RNAs and plays a crucial role in a number of regulatory
processes. Recently several MFE folding algorithms for predicting the
joint structure of two interacting RNA molecules have been proposed.
Here joint structure means that in a diagram representation the
intramolecular bonds of each partner are pseudoknot-free, that the
intermolecular binding pairs are noncrossing, and that there is no
so-called ``zig-zag'' configuration.
This paper presents the combinatorics of RNA interaction structures
including their generating function, singularity analysis as well as
explicit recurrence relations. In particular, our results imply simple
asymptotic formulas for the number of joint structures.
}

{\bf Keywords}: RNA-RNA interaction, Joint structure, Shape,
Symbolic enumeration, Singularity analysis.


\section{Introduction}


RNA-RNA binding is an important phenomenon observed in various classes
of non-coding RNAs and plays a crucial role in a number of regulatory
processes. Examples include the regulation of translation in both:
prokaryotes \citep{Vogel:07} and eukaryotes \citep{McManus,Banerjee},
the targeting of chemical modifications \citep{Bachellerie},
insertion editing \citep{Benne}, and transcriptional control
\citep{Kugel}.
More and more evidence suggests, that RNA-RNA interactions also
play a role for the functionality of long mRNA-like ncRNAs.
A common theme in many RNA classes, including miRNAs,
snRNAs, gRNAs, snoRNAs, and in particular many of the procaryotic small
RNAs, is the formation of RNA-RNA interaction structures that are much
more complex than simple complementary sense-antisense interactions.
The interaction between two RNAs is governed by the same physical
principles that determine RNA folding: the formation of specific
base pairs patterns whose energy is largely determined by base
pair stacking and loop strains.
Therefore, secondary structures are an appropriate level of description
to quantitatively understand the thermodynamics of RNA-RNA binding.

\begin{figure*}
\begin{tabular}{ccc}
\includegraphics[width=0.33\textwidth]{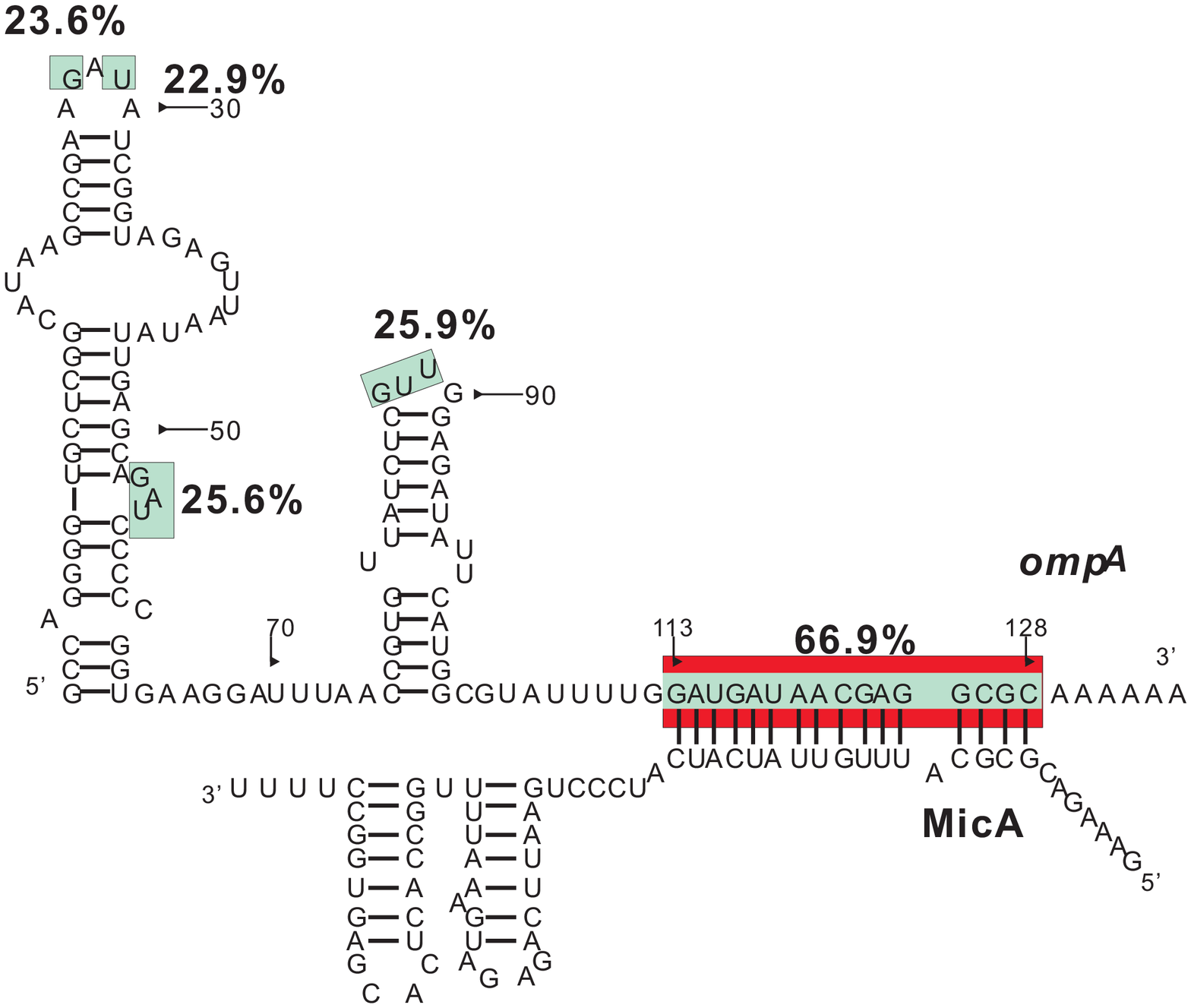}&
\includegraphics[width=0.3\textwidth]{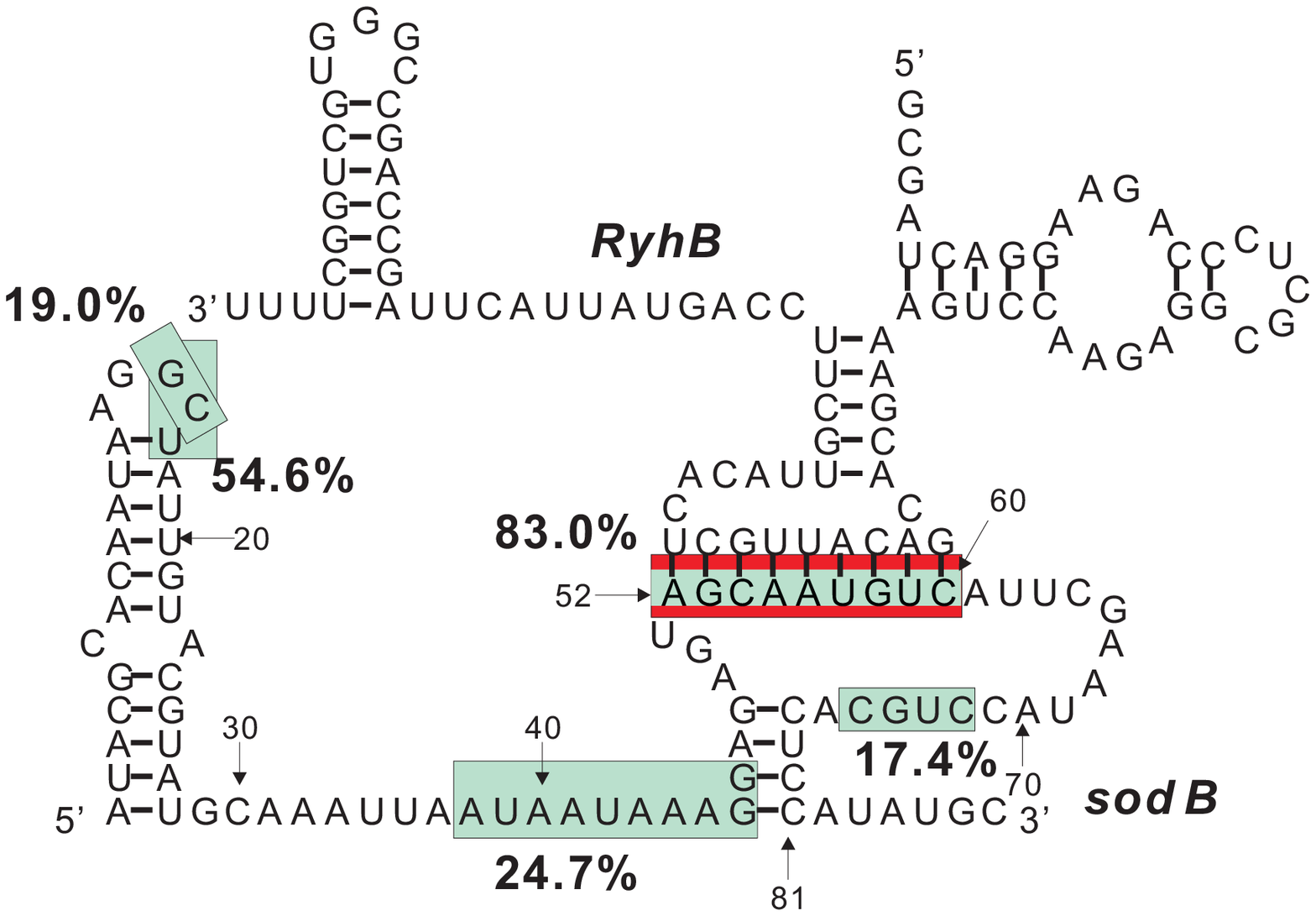}&
\includegraphics[width=0.3\textwidth]{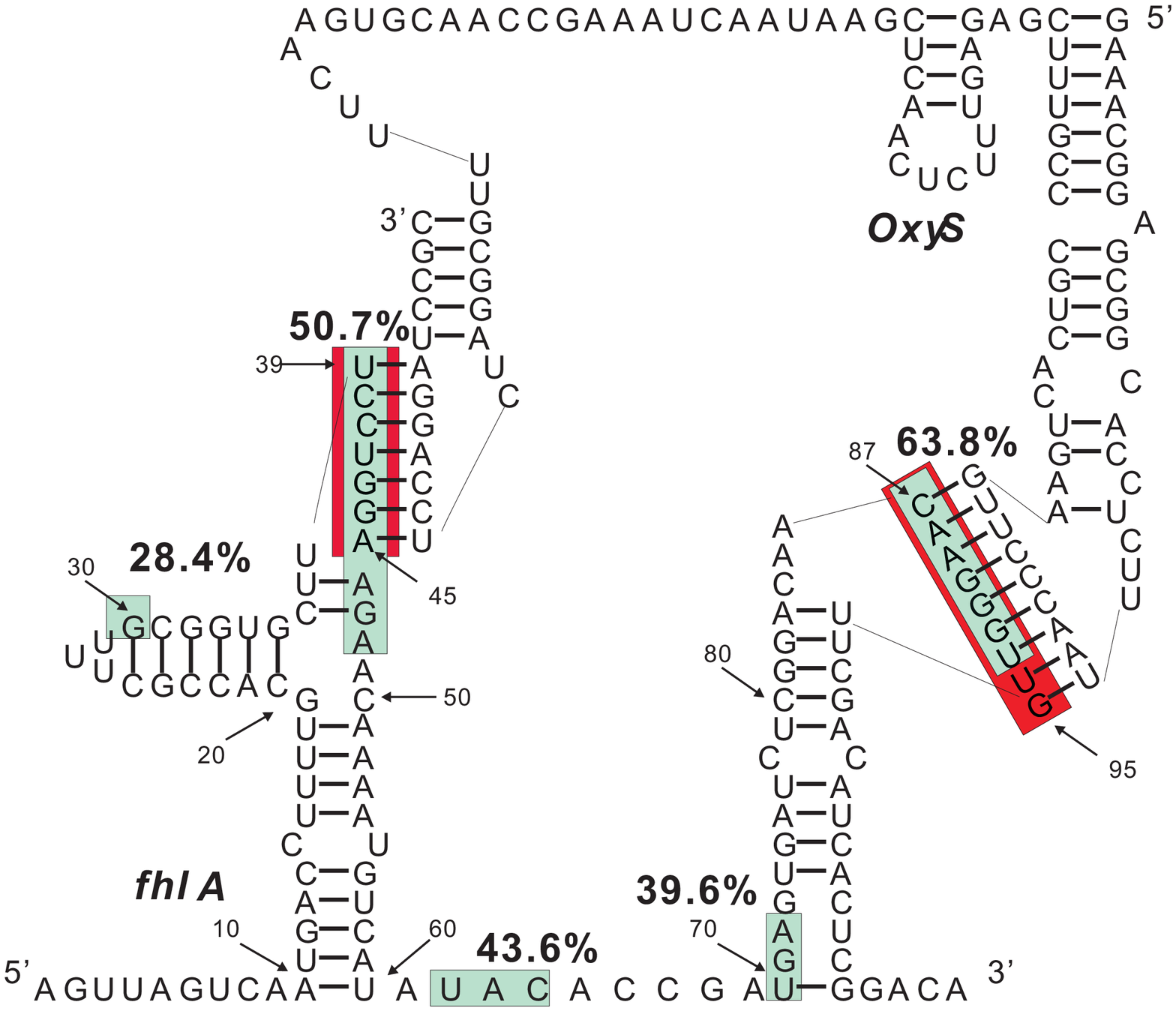}\\
\emph{ompA}-\emph{MicA} &
\emph{sodB}-\emph{RyhB} &
\emph{fhlA}-\emph{OxyS}\\
\end{tabular}
\caption{RNA-RNA interactions structures and their prediction.  The primary
  interaction region(s) are highlighted in red in the experimentally
  supported structural models from the literature:
  \emph{ompA}-\emph{MicA}: \citep{Udekwu:05};
  \emph{sodB}-\emph{RyhB}: \citep{Geissmann};
  \emph{fhlA}-\emph{OxyS}: \citep{Backofen}.
  Hybridization probabilities computed by \texttt{rip} are annotated
  by green boxes for regions with a probability larger than 10\%.
} \label{F:bac}
\end{figure*}

By restricting the space of allowed configurations, polynomial-time algorithms
on secondary structure level have been derived.
\citep{Pervouchine:04} and \citep{Alkan:06} proposed MFE folding
algorithms for predicting the \emph{joint structure} of two
interacting RNA molecules. In this model, ``joint structure'' means
that the intramolecular structures of each partner are
pseudoknot-free, that the intermolecular binding pairs are
noncrossing, and that there is no so-called ``zig-zag''
configuration, see Section~\ref{S:joint} for details.
This structure class seems to include
all major interaction complexes.
The optimal joint structure can be computed in $O(N^6)$ time and $O(N^4)$
space by means of dynamic programming \citep{Alkan:06,Pervouchine:04,rip:10,Backofen}.
More recently, extensions involving the partition function were proposed by
\citep{Backofen} (\texttt{piRNA}) and \citep{rip:09} (\texttt{rip}),
see Fig.~\ref{F:bac}.

In contrast to the situation for RNA secondary structures
\citep{Waterman:78,Waterman:94a}, little
is known about the joint structures that are the folding targets of \texttt{rip}
\citep{rip:10}.
This paper closes this gap and introduces the combinatorics of interaction
structures. We present the generating function of joint structures, its singularity
analysis as well as explicit recurrence relations. In particular, our results imply
simple formulas for the asymptotic number of joint structures.

The paper is organized as follows: in Section~\ref{S:facts} we provide several basic
fact and context. In Section~\ref{S:joint} we introduce joint structures along the
lines of \citep{rip:09}. In Section~\ref{S:shape} we follow the ideas of the paper
\citep{Modular} and consider shapes of joint structures.
In Section~\ref{S:inflate} we use shapes in order to compute the generating function
of joint structures and Section~\ref{S:asy} deals with the singularity analysis.
We then integrate our results in Section~\ref{S:discussion}. Finally we present
additional results in Section~\ref{S:later}.


\section{Some basic facts}\label{S:facts}


\subsection{Singularity analysis}

Let $f(z)=\sum_{n\geq 0}a_n\, z^n$ be a generating function with
nonnegative coefficients and a radius of convergence $R>0$. In light
of the fact that explicit formulas for the coefficients $a_n$ can be
very complicated or even impossible to obtain, we switch over to
investigate the estimation of $a_n$ in terms of the exponential
factor $\gamma$ and the subexponential factor $P(n)$, that is, $a_n
\sim P(n)\,\gamma^n$. The derivation of exponential growth rate and
subexponential factor is mainly based on singularity analysis.
Singularity analysis is a framework that allows to analyze the
asymptotics of these coefficients. The key to obtain the asymptotic
information about the coefficients of a generating function is its
dominant singularities, which raises the question on how to locate
them. In the particular case of power series with nonnegative
coefficients and a radius of convergence $R>0$, a theorem of
Pringsheim \citep{Flajolet:07a,Tichmarsh:39}, guarantees a positive real
dominant singularity at $z=R$. As we are dealing here with
combinatorial generating functions we always have this dominant
singularity. Furthermore for all our generating functions it is the
unique dominant singularity. The class of theorems that deal with
the deduction of information about coefficients from the generating
function are called transfer-theorems \citep{Flajolet:07a}.\\
To be precise, we say a function $f(z)$ is $\Delta_\rho$ analytic at
its dominant singularity $z=\rho$, if it analytic in some domain
$\Delta_\rho(\phi,r)=\{ z\mid \vert z\vert < r, z\neq \rho,\, \vert
{\rm Arg}(z-\rho)\vert >\phi\}$, for some $\phi,r$, where $r>|\rho|$
and $0<\phi<\frac{\pi}{2}$. We use the notation
\begin{equation*}
\left(f(z)=\Theta\left(g(z)\right) \
\text{\rm as $z\rightarrow \rho$}\right)\  \Longleftrightarrow \
\left(f(z)/g(z)\rightarrow c\ \text{\rm as $z\rightarrow \rho$}\right),
\end{equation*}
where $c$ is some constant. Let $[z^n]f(z)$ denote the coefficient
of $z^n$ in the power series expansion of $f(z)$ at $z=0$. Since the
Taylor coefficients have the property
\begin{equation*}
\forall \,\gamma\in\mathbb{C}\setminus 0;\quad [z^n]f(z)=\gamma^n
[z^n]f\left(\frac{z}{\gamma}\right),
\end{equation*}
We can, without loss of generality, reduce our analysis to the case
where $z=1$ is the unique dominant singularity. The next theorem
transfers the asymptotic expansion of a function around its unique
dominant singularity to the asymptotic of the function's
coefficients.
\begin{theorem}\label{T:transfer1}\citep{Flajolet:07a}
Let $f(z)$ be a $\Delta_1$ analytic function at its unique dominant
singularity $z=1$. Let
$$g(z)=(1-z)^{\alpha}\log^{\beta}\left(\frac{1}{1-z}\right)
,\quad\alpha,\beta\in \mathbb{R}.
$$
That is we have in the intersection of a neighborhood of $1$
\begin{equation}\label{E:transfer1}
f(z) = \Theta(g(z)) \quad \text{\it for } z\rightarrow 1.
\end{equation}
Then we have
\begin{equation}
[z^n]f(z)= \Theta\left([z^n]g(z)\right).
\end{equation}
\end{theorem}
\begin{theorem}\label{T:transfer2}\citep{Flajolet:07a}
Suppose $f(z)=(1-z)^{-\alpha}$, $\alpha\in\mathbb{C}\setminus
\mathbb{Z}_{\le 0}$, then
\begin{equation}
\begin{aligned}
[z^n]\, f(z)  \sim & \frac{n^{\alpha-1}}{\Gamma(\alpha)}\left[
1+\frac{\alpha(\alpha-1)}{2n}+\frac{\alpha(\alpha-1)(\alpha-2)(3\alpha-1)}
{24 n^2}+\right. \\
&\qquad  \quad
\left.\frac{\alpha^2(\alpha-1)^2(\alpha-2)(\alpha-3)}{48n^3}+
O\left(\frac{1}{n^4}\right)\right].
\end{aligned}
\end{equation}
\end{theorem}


\subsection{Symbolic Enumeration}

Symbolic enumeration \citep{Flajolet:07a} plays an important role in the following
computations. We first introduce the notion of a combinatorial class. Let
$\mathbf{z}=(z_1,\ldots,z_d)$ be a vector of $d$ formal variables
and $\mathbf{k}=(k_1,\ldots,k_d)$ be a vector of integers of the
same dimension. We use the simplified notation
\begin{equation*}
\mathbf{z}^{\mathbf{k}}\colon=z_1^{k_1} \cdots z_d^{k_d}.
\end{equation*}
\begin{definition}
A combinatorial class of $d$ dimension, or simply a class, is an
ordered pair $(\mathcal{A},w_{\mathcal{A}})$ where $\mathcal{A}$ is
a finite or denumerable set and a size-function
$w_{\mathcal{A}}\colon \mathcal{A}\longrightarrow \mathbb{Z}_{\geq
0}^d$ satisfies that $w_{\mathcal{A}}^{-1}(\mathbf{n})$ is finite
for any $\mathbf{n}\in\mathbb{Z}_{\geq 0}^d$.
\end{definition}
Given a class $(\mathcal{A},w_{\mathcal{A}})$, the size of an
element $a\in\mathcal{A}$ is denoted by $w_{\mathcal{A}}(a)$, or
simply $w(a)$. We consistently denote by $\mathcal{A}_{\mathbf{n}}$
the set of elements in $\mathcal{A}$ that have size $\mathbf{n}$ and
use the same group of letters for the cardinality
$A_{\mathbf{n}}=\vert \mathcal{A}_{\mathbf{n}}\vert$. The sequence
$\{A_{\mathbf{n}}\}$ is called the counting sequence of class
$\mathcal{A}$. The generating function of a class
$(\mathcal{A},w_{\mathcal{A}})$ is given by
\[
\mathbf{A}(\mathbf{z})=\sum_{a\in\mathcal{A}}\mathbf{z}^{w_{\mathcal{A}}(a)}=
\sum_{\mathbf{n}}A_{\mathbf{n}}\, \mathbf{z}^{\mathbf{n}}.
\]
There are two special classes: $\mathcal{E}$ and $\mathcal{Z}_i$
which contain only one element of size $\bf{0}$ and $\mathbf{e}_i$,
respectively. In particular, the generating functions of the classes
$\mathcal{E}$ and $\mathcal{Z}_i$ are
\[
{\bf{E}}(\mathbf{z})=1\quad\text{and}\quad {\bf{Z}}_i(\mathbf{z})=z_i.
\]
We adhere in the following to a systematic naming convention: classes, their
counting sequences, and their generating functions are
systematically denoted by the same groups of letters: for instance,
$\mathcal{C}$ for a class, $\{C_{\mathbf{n}}\}$ for the counting
sequence, and $\mathbf{C}(\mathbf{z})$ for its
generating function.
Let $\mathcal{A}$ and
$\mathcal{B}$ be combinatorial classes of $d$ dimension. Suppose
$\mathcal{A}_{i}$ are combinatorial classes of $1$ dimension. We
define
\begin{itemize}
\item $(\mathcal{A}_1,\mathcal{A}_2):=\{c=(a_1,a_2)\mid
a_i\in\mathcal{A}_i\}$ and for
$c=(a_1,a_2)\in(\mathcal{A}_1,\mathcal{A}_2)$
\[
w_{(\mathcal{A}_1,\mathcal{A}_2)}(c)=
(w_{\mathcal{A}_1}(a_1),w_{\mathcal{A}_2}(a_2))),
\]
\item $\mathcal{A}+\mathcal{B}:=\mathcal{A}\cup\mathcal{B}$, if
$\mathcal{A}\cap\mathcal{B}=\varnothing$ and for
$c\in\mathcal{A}+\mathcal{B}$,
\[
w_{\mathcal{A}+\mathcal{B}}(c)=\left\{
\begin{aligned}
&w_{\mathcal{A}}(c)\quad \text{if}\ c\in\mathcal{A}\\
&w_{\mathcal{B}}(c)\quad \text{if}\ c\in\mathcal{B},
\end{aligned} \right.
\]
\item $\mathcal{A}\times\mathcal{B}:=\{c=(a,b)\mid
a\in\mathcal{A},b\in\mathcal{B}\}$ and for
$c\in\mathcal{A}\times\mathcal{B}$,
\[
w_{\mathcal{A}\times\mathcal{B}}(c)=
w_{\mathcal{A}}(a)+w_{\mathcal{B}}(b),
\]
\item
$\textsc{Seq}(\mathcal{A}):={\mathcal{E}}+\mathcal{A}+(\mathcal{A}\times\mathcal{A})+
(\mathcal{A}\times\mathcal{A}\times\mathcal{A})+\cdots$.
\end{itemize}
Plainly, $\textsc{Seq}(\mathcal{A})$ defines a proper combinatorial
class if and only if $\mathcal{A}$ contains no element of size $0$.
We immediately observe
\begin{proposition}\label{P:Symbolic}
Suppose $\mathcal{A}$, $\mathcal{B}$ and $\mathcal{C}$ are
combinatorial classes of $d$ dimension having the generating functions
$\mathbf{A}(\mathbf{z})$, $\mathbf{B}(\mathbf{z})$ and
$\mathbf{C}(\mathbf{z})$. Let $\mathcal{A}_{i}$ be combinatorial
classes of $1$ dimension having the generating functions
$\mathbf{A}_i(z)$. Then\\
{\sf(a)} $\mathcal{C}=(\mathcal{A}_1,\mathcal{A}_2,\ldots,\mathcal{A}_d)
\Longrightarrow \mathbf{C}(\mathbf{z})=\mathbf{A}_1(z_1)\,
\mathbf{A}_2(z_2)\ldots \mathbf{A}_d(z_d)$\\
{\sf (b)} $\mathcal{C} =\mathcal{A}+\mathcal{B} \Longrightarrow
\mathbf{C}(\mathbf{z})=\mathbf{A}(\mathbf{z})+\mathbf{B}(\mathbf{z})$\\
{\sf (c)} $\mathcal{C}=\mathcal{A}\times\mathcal{B}\Longrightarrow
\mathbf{C}(\mathbf{z})=\mathbf{A}(\mathbf{z})\cdot\mathbf{B}(\mathbf{z})$\\
{\sf (d)} $\mathcal{C}=\textsc{Seq}(\mathcal{A})\Longrightarrow
\mathbf{C}(\mathbf{z})=\frac{1}{1-\mathbf{A}(\mathbf{z})}.$
\end{proposition}


\subsection{Secondary structures}


Let $f(n)$ denote the number of all noncrossing matchings of $n$
arcs having generating function ${\bf F}(z)= \sum f(n)\, z^n$.
Recursions for $f(n)$ allow us to derive
\begin{equation*}
z\, {\bf F}(z)^2 -{\bf F}(z) + 1 = 0,
\end{equation*}
that is we have
\begin{equation*}
{\bf F}(z) = \frac{1-\sqrt{1-4z}}{2z}.
\end{equation*}
Let $\mathcal{T}_{\sigma}$ denote the combinatorial
class of $\sigma$-canonical secondary structures having arc-length
$\geq 2$ and $T_{\sigma}(n)$ denote the number of
all $\sigma$-canonical secondary structures with $n$ vertices having
arc-length $\geq 2$ and
\begin{equation*}
{\bf T}_{\sigma}(z)= \sum T_{\sigma}(n)\, z^n.
\end{equation*}
\begin{theorem}\label{T:SecondSig}
Suppose $\sigma\in \mathbb{N}$, $\sigma\geq 1$ and $u_\sigma(z)
=\frac{(z^2)^{\sigma-1}}{z^{2\sigma}-z^2+1}$. Then we have
\begin{eqnarray*}
{\bf T}_{\sigma}(z) & = & \frac{1}
{u_\sigma(z) z^2-z+1}{\bf F}\left(\left(\frac{\sqrt{u_\sigma(z)}z}
{\left(u_\sigma(z) z^{2}-z+1\right)}\right)^2\right).
\end{eqnarray*}
where
\begin{equation*}
{\bf F}(z)   = \frac{1-\sqrt{1-4z}}{2z}.
\end{equation*}
\end{theorem}
Since ${\bf F}(z)$ is algebraic and $u_\sigma(z)$ is a rational function,
Theorem~\ref{T:SecondSig} implies that ${\bf T}_{\sigma}(z)$ is an algebraic
function for any $\sigma$.


\section{Joint Structures }\label{S:joint}

\begin{figure}
\includegraphics[width=1\textwidth]{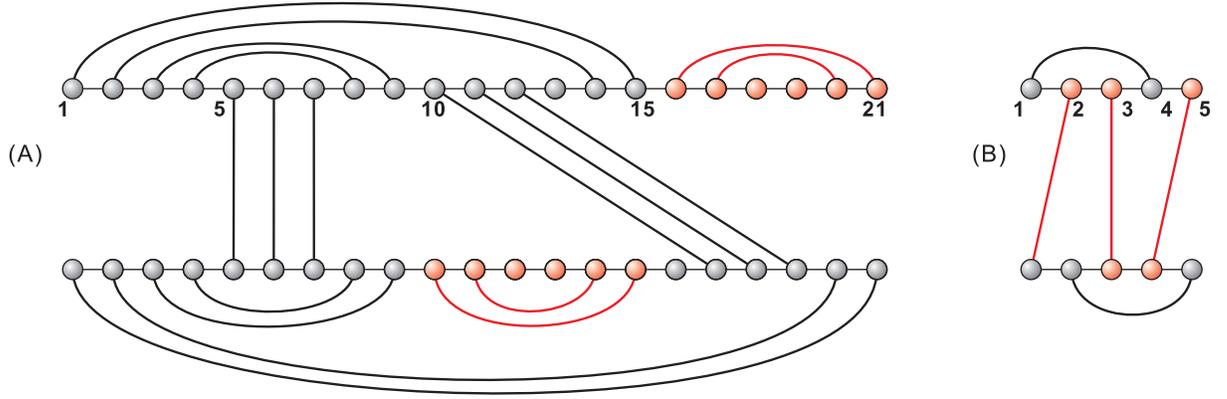}
\caption{\textsf{(A)}: The joint structure $J(R,S,I)$ with arc-length
$\geq 3$, interior stack-length $\geq 2$, exterior stack-length $\geq 3$.
Secondary segments (red): the subgraphs $R[16,21]$ and $S[10,15]$.
Ancestors and descendants: for the exterior arc $R_5 S_5$, we have the
following sets of $R$-ancestors and $S$-ancestors of $R_5 S_5$: $\{ R_1 R_{15},
R_2 R_{14}, R_3 R_{9}, R_4 R_{8},\}$ and $\{ S_1 S_{21}, S_2 S_{20}, S_3 S_{9},
S_4 S_{8},\}$. The exterior arc $R_5 S_5$ is a common descendant of $R_1 R_{15}$
and $S_3 S_{9}$, while $R_{10} S_{17}$ is not.
Subsumed arcs: $R_1 R_{15}$ subsumes $S_3 S_9$ and $S_1 S_{21}$.
\textsf{(B)}: A zigzag, generated by $R_2 S_1$, $R_3 S_3$ and $R_5 S_4$.
}
\label{F:JSexample}
\end{figure}

Given two RNA sequences $R$ and $S$ with $n$ and $m$ vertices, we
index the vertices such that $R_1$ is the $5'$ end of $R$ and $S_1$
is the $3'$ end of $S$. We refer to the $i$th vertex in $R$ by $R_i$
and the subgraph induced by $\{R_i,\ldots,R_j\}$ by $R[i,j]$. The
intramolecular base pair can be represented by an arc (interior),
with its two endpoints contained in either $R$ or $S$. Similarly,
the extramolecular base pair can be represented by an arc (exterior)
with one of its endpoints contained in $R$ and the other in $S$. A
pre-structure, $G(R,S,I)$, is a graph consisting of two secondary
structures $R$ and $S$ with a set $I$ of noncrossing exterior arcs.
When representing arc-configurations, we draw all $R$-arcs in the
upper-halfplane and all $S$-arcs in the lower-halfplane, see
Fig.~\ref{F:JSexample},~{\bf (A)}.

The subgraph $R[i,j]$ ($S[i',j']$) is called secondary segment if
there is no exterior arc $R_k S_{k'}$ such that $i \leq k \leq j$
($i' \leq k' \leq j'$), see Fig.~\ref{F:JSexample}, {\bf (A)}. An
interior arc $R_i R_j$ is an $R$-ancestor of the exterior arc $R_k
S_{k'}$ if $i<k<j$. Analogously, $S_{i'} S_{j'}$ is an $S$-ancestor
of $R_k S_{k'}$ if $i'<k'<j'$. We also refer to $R_k S_{k'}$ as a
descendant of $R_i R_j$ and $S_{i'} S_{j'}$ in this situation, see
Fig.~\ref{F:JSexample}, {\bf (A)}. Furthermore, we call $R_i R_j$
and $S_{i'} S_{j'}$ dependent if they have a common descendant and
independent,otherwise. Let $R_i R_j$ and $S_{i'} S_{j'}$ be two
dependent interior arcs. Then $R_i R_j$ subsumes $S_{i'} S_{j'}$, or
$S_{i'} S_{j'}$ is subsumed in $R_i R_j$, if for any $R_k S_{k'}\in
I$, $i'<k'<j'$ implies $i<k<j$, that is, the set of descendants of
$S_{i'} S_{j'}$ is contained in the set of descendants of $R_i R_j$,
see Fig.~\ref{F:JSexample},~{\bf (A)}. A zigzag is a subgraph
containing two dependent interior arcs $R_{i_1} R_{j_1}$ and
$S_{i_2} S_{j_2}$ neither one subsuming the other, see
Fig.~\ref{F:JSexample}, {\bf (B)}. A joint structure $J(R,S,I)$ is a
zigzag-free pre-structure, see Fig.~\ref{F:JSexample},~{\bf (A)}.\\
We denote the combinatorial class of all joint structures by
$\mathcal{J}$. We can define the size-function as follows:
$w_{\mathcal{J}}(J(R,S,I))= (n,m,h)$, where $n$ and $m$ denote the
number of vertices in the top and bottom sequence and $h$ denotes
the number of exterior arcs in the joint structure. We denote by
$\mathcal{J}(n,m,h)$ the subset of $\mathcal{J}$ which contains all
the joint structures of the size $(n,m,h)$ and set the counting
sequence $J(n,m,h)=\vert \mathcal{J}(n,m,h) \vert$. The generating
function of the class $\mathcal{J}$ is given by
\[
\mathbf{J}(x,y,z)=\sum_{n,m}J(n,m,h)\, x^n y^m z^h.
\]
We next specify some notation
\begin{itemize}
\item an interior arc (or simply arc) of length $\lambda$ is an arc $R[i,j]$
    ($S[i',j']$) where $j-i=\lambda$ ($j'-i'=\lambda$),
\item an interior stack (or simply stack) of length $\sigma$ is a maximal
      sequence of ``parallel'' interior arcs,
    \begin{eqnarray*}
    (R_i R_j, R_{i+1} R_{j-1}, \ldots , R_{i+\sigma-1} R_{j-\sigma+1} ) &
\text{\rm or}& \\
    (S_i S_j, S_{i+1} S_{j-1}, \ldots , S_{i+\sigma-1} S_{j-\sigma+1} ), &&
    \end{eqnarray*}
\item an exterior stack of length $\tau$ is a maximal sequence of ``parallel''
     exterior arcs,
    \begin{equation*}
    (R_i S_{i'}, R_{i+1} S_{i'+1}, \ldots , R_{i+\tau-1} S_{i'+\tau-1} ).
    \end{equation*}
\end{itemize}
Let $\mathcal{J}_{\sigma,\tau}^{[\lambda]}$ denote the class of all
joint structures with arc-length $\geq \lambda$, interior
stack-length $\geq \sigma$, exterior stack-length $\geq \tau$.
Similarly, we can define its counting sequence
$J_{\sigma,\tau}^{[\lambda]}(n,m,h)$ and generating function
$\mathbf{J}_{\sigma,\tau}^{[\lambda]}(x,y,z)$. In case of $\lambda =2$, we
omit $\lambda$ in the notation. If there is no restriction on the
interior and exterior stack-length, we also omit further indices.
In the particular case $\sigma = \tau$, we just
write $\sigma$ in the notation and omit $\tau$. In
Fig.~\ref{F:JSexample}, {\bf (A)}, we give an example of joint
structure with arc-length $\geq 3$, interior stack-length $\geq 2$ and
exterior stack-length $\geq 3$.\\

We denote the subgraph of a joint structure $J(R,S,I)$ induced by a pair of
subsequences $\{R_i,R_{i+1},\ldots,R_j\}$ and $\{S_{i'},S_{i'+1},\ldots,
S_{j'}\}$ the block $J_{i,j;i',j'}$.
Given a joint structure $J(R,S,I)$, a tight structure of $J(R,S,I)$ is
the minimal block $J_{i,j;i',j'}$ containing all the $R$-ancestors and
$S$-ancestors of any exterior arc in $J_{i,j;i',j'}$ and all the descendants
of any interior arc in $J_{i,j;i',j'}$. In the following, a tight
structure is denoted by $J^{T}_{i,j;i',j'}$. In particular, we
denote the joint structure $J(R,S,I)$ by $J^T (R,S,I)$ if $J(R,S,I)$
is a tight structure of itself. For any joint structure, there are
only four types of tight structures $J^T_{i,j;i',j'}$, that is
$\{\circ,\triangledown,\vartriangle,\square \}$, denoted by $J^{\{
\circ,\triangledown,\vartriangle,\square
 \}}_{i,j;i',j'}$, respectively. The four types of tight structures $J^{\{
\circ,\triangledown,\vartriangle,\square
 \}}_{i,j;i',j'}$ are defined as follows:
\begin{align*}
\circ           &\colon \{R_i S_{i'}\}=J^{\circ}_{i,j;i',j'}
\quad\text{and}\quad i=j\ ,\ i'=j';\\
\triangledown   &\colon R_i R_j \in J^{\triangledown}_{i,j;i',j'}
\quad\text{and}\quad S_{i'} S_{j'} \notin J^{\triangledown}_{i,j;i',j'};\\
\vartriangle    &\colon S_{i'} S_{j'} \in J^{\vartriangle}_{i,j;i',j'}
\quad\text{and}\quad R_i R_j \notin J^{\vartriangle}_{i,j;i',j'};\\
\square         &\colon \{R_i R_j,S_{i'} S_{j'}\} \in J^{\square}_{i,j;i',j'}.
\end{align*}

The key function of tight structures is that they are the building blocks for
the decomposition of joint structures.
\begin{proposition}\label{P:Decomposition}\citep{rip:09}
Let $J(R,S,I)$ be a joint structure. Then
\begin{enumerate}
\item any exterior arc $R_k S_{k'}$ in $J(R,S,I)$ is contained in a
      unique tight structure.
\item $J(R,S,I)$ decomposes into a unique collection of tight structures
      and maximal secondary segments.
\end{enumerate}
\end{proposition}


\section{Shapes}\label{S:shape}

\begin{figure}
\includegraphics[width=1\textwidth]{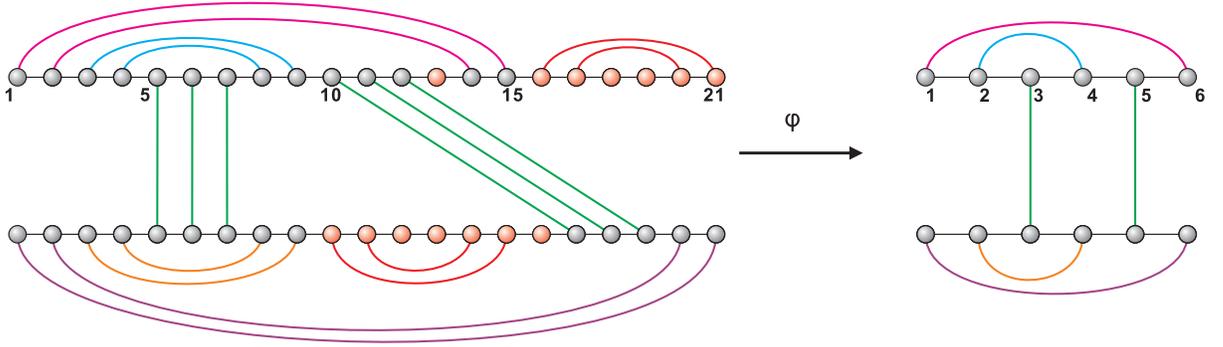}
\caption{Joint structures and their shapes: a joint structure (left) is
projected into its shape (right). } \label{F:ShapeExample}
\end{figure}

\begin{definition}{\bf (Shape)}
A shape is a joint structure containing no secondary segments in
which each interior stack and each exterior stack have length
exactly one.
\end{definition}
Let $\mathcal{G}$ denote the combinatorial class of shapes. Given a joint
structure, we can obtain its shape by first removing all secondary segments
and second collapsing any stacks into a single arc.
That is, we have a map $\varphi \colon {\mathcal J}\rightarrow
{\mathcal G}$, see Fig.~\ref{F:ShapeExample}.
Let $G(t_1,t_2,h)$ denote the number of shapes having $t_1$ arcs in
the top sequence, $t_2$ arcs in the bottom and $h$ exterior arcs
having the generating function
\begin{equation*}
{\bf G}(u,v,z)=\sum G(t_1,t_2,h)\, u^{t_1} v^{t_2} z^h.
\end{equation*}
We next introduce tight shapes, double tight shapes,
interaction segments, closed shapes and right closed shapes:
\begin{itemize}
\item A tight shape is tight as a structure.
    Let $\mathcal{G}^T$ denote the class of tight shapes by and
    $G^T (t_1,t_2,h)$ denote the number of tight shapes having
    $t_1$ arcs in the top sequence, $t_2$ arcs in the bottom
    and $h$ exterior arcs having the generating function
    \begin{equation*}
    {\bf G}^T(u,v,z)=\sum G^T(t_1,t_2,h)\, u^{t_1} v^{t_2} z^h.
    \end{equation*}
    Any tight shape, comes as exactly one of the four types
    $\{\circ,\triangledown,\vartriangle,\square\}$.
    The corresponding classes and generating functions are
    defined accordingly,
    $\mathcal{G}^{\{\circ,\triangledown,\vartriangle,\square\}}$
    and ${\bf
    G}^{\{\circ,\triangledown,\vartriangle,\square\}}(x,y,z)$
    respectively,
\item A double tight shape is a shape whose leftmost and rightmost
    blocks are tight structures. Let $\mathcal{G}^{DT}$ denote
    the class of double tight shapes by and $G^{DT} (t_1,t_2,h)$
    denote the number of double tight shapes having $t_1$ arcs in the
    top sequence, $t_2$ arcs in the bottom and $h$ exterior
    arcs having the generating function
    \begin{equation*}
    {\bf G}^{DT}(u,v,z)=\sum G^{DT}(t_1,t_2,h)\, u^{t_1} v^{t_2} z^h,
    \end{equation*}
\item A closed shape is a tight shape of type $\{\triangledown,
    \vartriangle,\square\}$. Let $\mathcal{G}^C$ denote the class of
    closed shapes and $G^C (t_1,t_2,h)$ denote the number of closed
    shapes having $t_1$ arcs in the top sequence, $t_2$ arcs in the
    bottom and $h$ exterior arcs having the generating function
    \begin{equation*}
    {\bf G}^C(u,v,z)=\sum G^C(t_1,t_2,h)\, u^{t_1} v^{t_2} z^h,
    \end{equation*}
\item A right closed shape is a shape whose rightmost block is
    a closed shape rather than an exterior arc. Let
    $\mathcal{G}^{RC}$ denote the class of right close shapes
    and $G^{RC} (t_1,t_2,h)$ denote the number of right close
    shapes having $t_1$ arcs in the top sequence, $t_2$ arcs
    in the bottom and $h$ exterior arcs having the generating function
    \begin{equation*}
    {\bf G}^{RC}(u,v,z)=\sum G^{RC}(t_1,t_2,h)\, u^{t_1} v^{t_2} z^h,
    \end{equation*}
\item In a shape, an interaction segment is an empty structure or an tight
    structure of type $\circ$ (an exterior arc). We denote the
    class of interaction segment by $\mathcal{I}$ and the associated
    generating function by ${\bf I}(x,y,z)$. Obviously, ${\bf
    I}(x,y,z)=1+z$.
\end{itemize}

\begin{figure}
\includegraphics[width=0.5\textwidth]{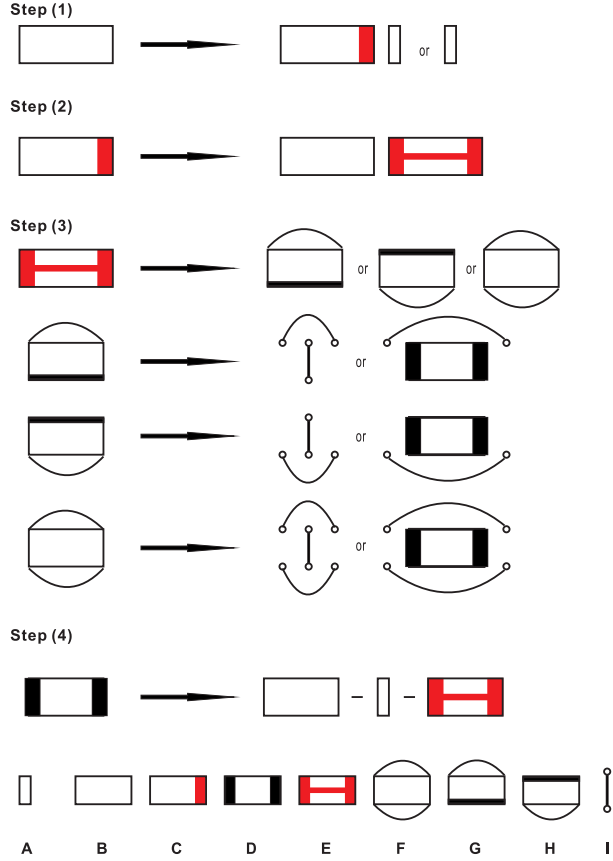}
\caption{The shape-grammar. The notations of structural components are
explained in the panel below.
{\bf{A}}: interaction segment;
{\bf{B}}: arbitrary shape $G(R,S,I)$; {\bf{C}}: right close shape
$G^{RC}(R,S,I)$; {\bf{D}}: double tight shape $G^{DT}(R,S,I)$;
{\bf{E}}: close shape $G^{C}(R,S,I)$; {\bf{F}}: type $\square$
tight shape $G^{\square}(R,S,I)$; {\bf{G}}: type $\triangledown
$ tight shape $G^{\triangledown }(R,S,I)$; {\bf{H}}: type
$\vartriangle$ tight shape $G^{\vartriangle}(R,S,I)$; {\bf{I}}:
type $\circ$ tight shape $G^{\circ}(R,S,I)$; }
\label{F:ShapeDecomGram}
\end{figure}

\begin{theorem}\label{T:Shape}
The generating function ${\bf G}(u,v,z)$ of shapes satisfies
\begin{equation}\label{E:Shape}
{\bf A}(u,v,z){\bf G}(u,v,z)^2+ {\bf B}(u,v,z){\bf G}(u,v,z)+
{\bf C}(u,v,z)=0,
\end{equation}
where
\begin{equation}\label{E:erni}
\begin{aligned}
{\bf A}(u,v,z)  &= (u+v+uv)(z+1),\\
{\bf B}(u,v,z)  &= -((u+v+uv)(z+2)+1),\\
{\bf C}(u,v,z)  &= (1+u)(1+v)(1+z) .
\end{aligned}
\end{equation}
\end{theorem}
\begin{proof}
Proposition~\ref{P:Decomposition} implies that any shape can be
decomposed into a unique collection of tight shapes. Furthermore,
each shape can be decomposed into a unique collection of close shapes
and exterior arcs. We decompose a shape in four steps, see
Fig.~\ref{F:ShapeDecomGram}. We translate each decomposition step into
the construction of combinatorial classes in the language of symbolic
enumeration. \\
{\bf Step (1):} we decompose a shape into a right closed shape and
rightmost interaction segment. We generate
$\mathcal{G}=\mathcal{G}^{RC}\times \mathcal{I}+ \mathcal{I}$.\\
It follows from Proposition~\ref{P:Symbolic} that
\begin{equation}\label{E:Shape-1}
{\bf G}(x,y,z)={\bf G}^{RC}(x,y,z)\cdot {\bf I}(x,y,z)+{\bf I}(x,y,z).
\end{equation}
{\bf Step (2):} we decompose a right closed shape into the rightmost
closed shape and the rest, deriving
\[
\mathcal{G}^{RC}=\mathcal{G}\times \mathcal{G}^{C},
\]
whence
\begin{equation}\label{E:Shape-2}
{\bf G}^{RC}(x,y,z)={\bf G}(x,y,z)\cdot {\bf G}^{C}(x,y,z).
\end{equation}
{\bf Step (3):} we decompose a closed shape depending on its type. The
decomposition operation in this step can be viewed as the
"removal" of an interior arc. We derive
\begin{eqnarray*}
\mathcal{G}^{C}               &=&
\mathcal{G}^{\triangledown}+\mathcal{G}^{\vartriangle}
+\mathcal{G}^{\square}\\
\mathcal{G}^{\triangledown}   &=&
(\mathcal{Z},\mathcal{E},\mathcal{Z}) +
(\mathcal{Z},\mathcal{E},\mathcal{E}) \times \mathcal{G}^{DT}\\
\mathcal{G}^{\vartriangle}   &=&
(\mathcal{E},\mathcal{Z},\mathcal{Z})+
(\mathcal{E},\mathcal{Z},\mathcal{E}) \times \mathcal{G}^{DT}\\
\mathcal{G}^{\square}       &=&
(\mathcal{Z},\mathcal{Z},\mathcal{Z})+
(\mathcal{Z},\mathcal{Z},\mathcal{E}) \times \mathcal{G}^{DT}
\end{eqnarray*}
and obtain the generating functions
\begin{equation}\label{E:Shape-3}
\begin{aligned}
{\bf G}^{C}(x,y,z)                  &=  {\bf G}^{\triangledown}(x,y,z)+
{\bf G}^{\vartriangle}(x,y,z)+{\bf G}^{\square}(x,y,z)\\
{\bf G}^{\triangledown}(x,y,z)      &=  x\, z +  x \, {\bf G}^{DT}(x,y,z)\\
{\bf G}^{\vartriangle}(x,y,z)       &=  y\, z+  y \, {\bf G}^{DT}(x,y,z)\\
{\bf G}^{\square}(x,y,z)            &=  x\, y\, z+  x\, y \, {\bf G}^{DT}(x,y,z).
\end{aligned}
\end{equation}
{\bf Step (4):} the class of double tight shapes arising from Step
(3) can be obtained by excluding the class of interaction segment
and the class of closed shapes from the class of shapes. Similarly,
we have
\[
\mathcal{G}^{DT}=\mathcal{G}-\mathcal{I}-\mathcal{G}^C.
\]
The corresponding generating function accordingly satisfies
\begin{equation}\label{E:Shape-4}
{\bf G}^{DT}(x,y,z)={\bf G}(x,y,z) - {\bf I}(x,y,z) -{\bf G}^{C}(x,y,z).
\end{equation}
We proceed by solving the set of equations~(\ref{E:Shape-1})--(\ref{E:Shape-4}),
thereby deriving the functional equation eq.~(\ref{E:erni}) for
${\bf G}(x,y,z)$ and the theorem follows.
\end{proof}


\section{The generating function}\label{S:inflate}


We proceed by generating joint structures from shapes via inflation.
Let ${\mathcal J}_{\sigma,\tau}$ denote the class of joint
structures with arc-length $\geq 2$, interior stack-length $\geq
\sigma$, exterior stack-length $\geq \tau$. Let
$J_{\sigma,\tau}(n,m,h)$ denote the number of joint structures
in ${\mathcal J}_{\sigma,\tau}$ having $n$ vertices in the top,
$m$ vertices in the bottom and $h$ exterior arcs having the generating
function
\begin{equation*}
{\bf J}_{\sigma,\tau}(x,y,z) =\sum J_{\sigma,\tau}(n,m,h)\, x^{n} y^{m} z^h.
\end{equation*}
\begin{theorem}\label{T:RelateJG}
For $\sigma \geq 1,\tau \geq 1$ , we have
\begin{equation}
{\bf J}_{\sigma,\tau}(x,y,z) = {\bf T}_{\sigma}(x)\,{\bf T}_{\sigma}(y)
\,{\bf G}\!\big(\eta(x),\eta(y),\eta_0\big),
\end{equation}
where
\begin{eqnarray*}
\eta(w)     &=&
\frac{w^{2\sigma}\,{\bf T}_{\sigma}(w)^2}{1-w^2-w^{2\sigma}
({\bf T}_{\sigma}(w)^2-1)},\\
\eta_0      &=&  \frac{(x y z)^{\tau}{\bf T}_{\sigma}(x)\,
{\bf T}_{\sigma}(y)}{1- x y z -(x y z)^{\tau}({\bf T}_{\sigma}(x)
{\bf T}_{\sigma}(y)-1)}.
\end{eqnarray*}
\end{theorem}
\begin{proof}
Let ${\mathcal G}(t_1,t_2,h)$ denote the class of shapes having $t_1$
interior arcs in the top, $t_2$ interior arcs in the bottom and $h$ exterior
arcs. For any joint structure, we can obtain a unique shape in
$\mathcal{G}$ as follows:
\begin{enumerate}
\item Remove all secondary segments.
\item Contract each interior stack into one interior arc and each exterior
stack into one exterior arc.
\end{enumerate}
Then we have the surjective map
\begin{equation*}
\varphi\colon {\mathcal J}_{\sigma,\tau}\rightarrow {\mathcal G}.
\end{equation*}
Indeed, for any shape $\gamma$ in $\mathcal{G}$, we can construct
joint structures with arc-length $\geq 2$, stack-length $\geq
\sigma$, exterior stack-length $\geq \tau$. $\varphi\colon {\mathcal
J}_{\sigma,\tau}\rightarrow {\mathcal G}$, induces the partition
${\mathcal J}_{\sigma,\tau}=\dot\cup_\gamma\varphi^{-1}(\gamma)$.
Then we have
\begin{equation}
{\bf J}_{\sigma,\tau}(x,y,z) = \sum_{\gamma\in\,{\mathcal G}}
\mathbf{J}_\gamma(x,y,z).
\end{equation}
We proceed by computing the generating function
$\mathbf{J}_\gamma(x,y,z)$. We will construct
$\mathbf{J}_\gamma(x,y,z)$ via simpler combinatorial classes as
building blocks considering $\mathcal{M}_{\sigma}$
(stems), $\mathcal{K}_{\sigma}$ (stacks), $\mathcal{N}_{\sigma}$
(induced stacks), $\mathcal{R}$ (interior arcs) and
$\mathcal{T}_{\sigma}$ (secondary segments).
We inflate a shape $\gamma\in {\mathcal G}(t_1,t_2,h)$ to a joint structure in
three steps.\\
{\bf Step I:} we inflate any interior arc in $\gamma$ to a stack of
size at least $\sigma$ and subsequently add additional stacks. The
latter are called induced stacks and have to be separated by means
of inserting secondary segments, see Fig.~\ref{F:RelateJG-1}.
\begin{figure}
\includegraphics[width=1\textwidth]{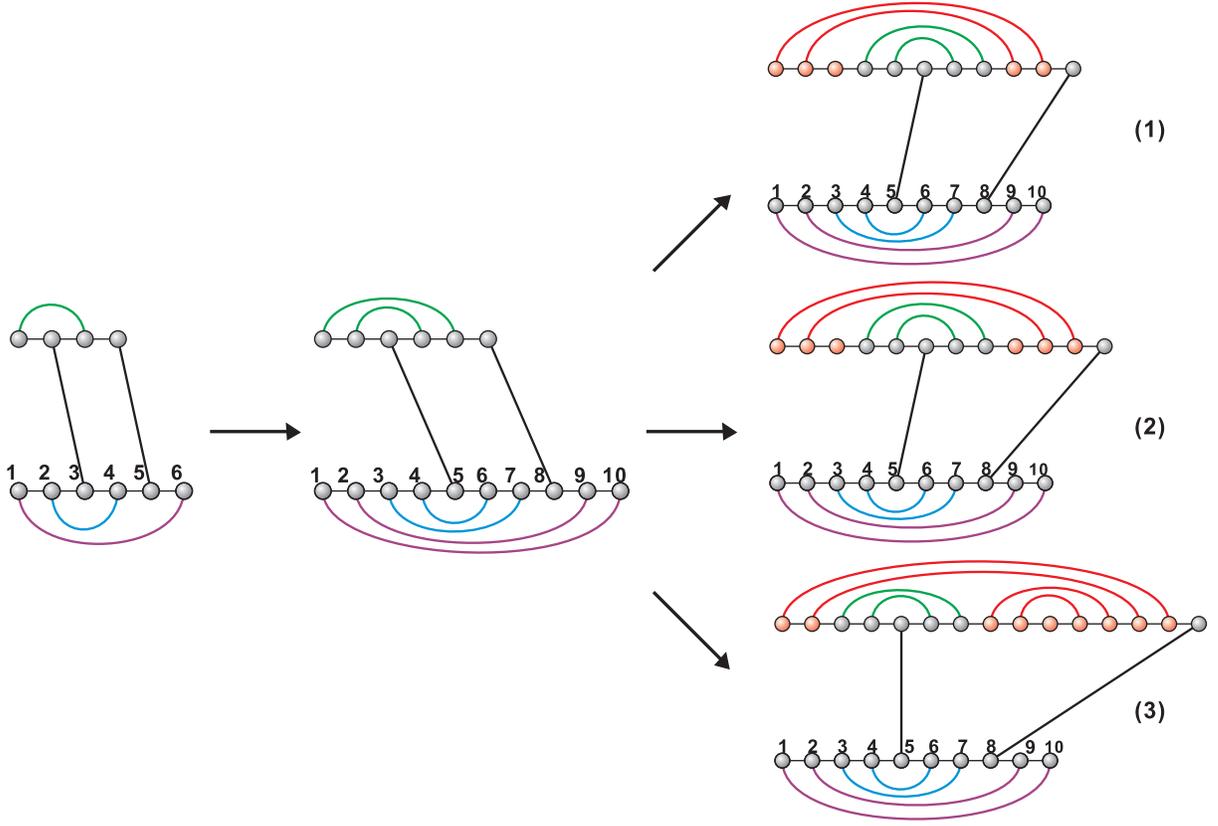}
\caption{Step I: a shape (left) is inflated to a joint structure with arc-length
$\geq 2$ and interior stack-length $\geq 2$. First, each interior arc in the shape
is inflated to a stack of size at least two (middle). Then the shape is inflated
to a new joint structure with arc-length $\geq 2$ and interior stack-length
$\geq 2$ (right) by adding one stack of size two. Note that there are three ways
to insert the secondary segments to separate the induced stacks (red).}
\label{F:RelateJG-1}
\end{figure}
Note that during this first inflation step no secondary segments,
other than those necessary for separating the nested stacks are
inserted. We generate
\begin{itemize}
\item secondary segments $\mathcal{T}_{\sigma}$ having stack-length $\geq \sigma$
having the generating function ${\bf T}_{\sigma}(z)$,
\item interior arcs $\mathcal{R}$ with generating function ${\bf R}(z)=z^2$,
\item stacks, i.e.~pairs consisting of the minimal sequence of arcs
$\mathcal{R}^\sigma$ and an arbitrary extension consisting of
arcs of arbitrary finite length
\begin{equation*}
\mathcal{K}_{\sigma}=\mathcal{R}^{\sigma}\times\textsc{Seq}\left(\mathcal{R}\right)
\end{equation*}
having the generating function
\begin{eqnarray*}
\mathbf{K}_\sigma(z) & = & z^{2\sigma}\cdot \frac{1}{1-z^2},
\end{eqnarray*}
\item induced stacks, i.e.~stacks together with at least one secondary segment on
either or both of its sides,
\begin{equation*}
\mathcal{N}_{\sigma}=\mathcal{K}_{\sigma}\times\left(
\mathcal{T}_{\sigma}^2 -1 \right),
\end{equation*}
having the generating function
\begin{equation*}
\mathbf{N}_\sigma(z)=\frac{z^{2\sigma}}{1-z^2}\left( {\bf T}_{\sigma}(z)^2-1\right),
\end{equation*}
\item stems, that is pairs consisting of stacks $\mathcal{K}_\sigma$
and an arbitrarily long sequence of induced stacks
\begin{equation*}
\mathcal{M}_\sigma=\mathcal{K}_{\sigma} \times
\textsc{Seq}\left(\mathcal{N}_{\sigma}\right),
\end{equation*}
having the generating function
\begin{eqnarray*}
\mathbf{M}_\sigma(z)=\frac{\mathbf{K}_\sigma(z)}{1-\mathbf{N}_\sigma(z)}=
\frac{\frac{z^{2\sigma}}{1-z^2}}
{1-\frac{z^{2\sigma}}{1-z^2}\left( {\bf T}_{\sigma}(z)^2-1 \right)}.
\end{eqnarray*}
\end{itemize}
Note that we inflate both: top as well as bottom sequences. The
corresponding generating function is
\begin{eqnarray*}
{\bf M}_{\sigma}(x)^{t_1}\,{\bf M}_{\sigma}(y)^{t_2}.
\end{eqnarray*}

{\bf Step II:} we inflate any exterior arc in $\gamma$ to an
exterior stack of size at least $\tau$ and subsequently add
additional exterior stacks. The latter are called induced exterior
stacks and have to be separated by means of inserting secondary
segments, see Fig.~\ref{F:RelateJG-2}.
\begin{figure}
\includegraphics[width=1\textwidth]{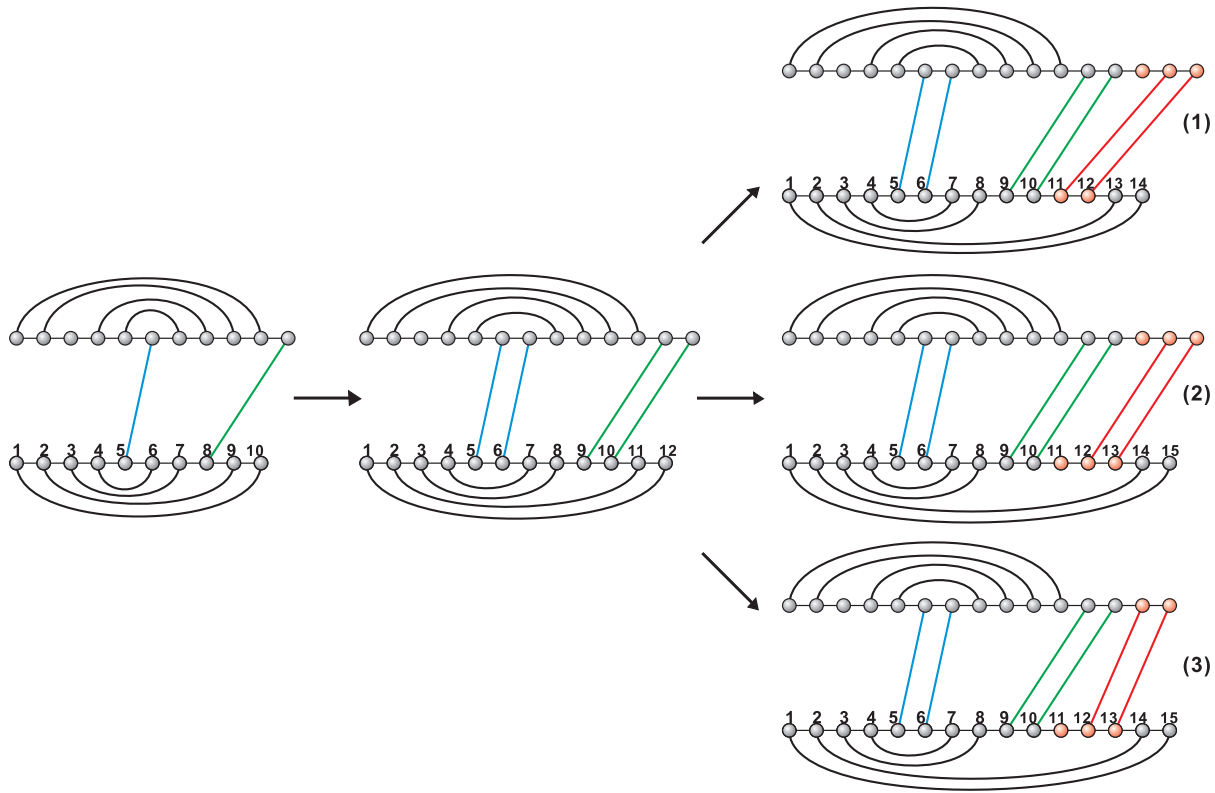}
\caption{Step II: a joint structure (left) obtained in {\sf (1)} in
Fig.~\ref{F:RelateJG-1} is inflated to a joint structure in ${\mathcal J}_{2,2}$.
First, each exterior arc in the joint structure is inflated to
an exterior stack of size at least two (middle), and then the structure is inflated
to a new joint structure in ${\mathcal J}_{2,2}$ (right) by adding
one exterior stack of size two. There are three ways to insert the secondary
segments to separate the induced exterior stacks (red).
} \label{F:RelateJG-2}
\end{figure}
Note that during this exterior-arc inflation step no secondary
segments, other than those necessary for separating the stacks are
inserted. We generate
\begin{itemize}
\item exterior arc $\mathcal{R}_0$ having the generating function
\begin{equation*}
{\bf R}_0= x y z,
\end{equation*}
\item exterior stacks, i.e.~pairs consisting of the minimal sequence of exterior arcs
$\mathcal{R}_0^\tau$ and an arbitrary extension consisting of
exterior arcs of arbitrary finite length
\begin{equation*}
\mathcal{K}'_{\tau}=\mathcal{R}_0^{\tau}\times\textsc{Seq}\left(\mathcal{R}_0\right)
\end{equation*}
having the generating function
\begin{eqnarray*}
\mathbf{K}'_{\tau} & = & (x y z)^{\tau}\cdot \frac{1}{1-x y z},
\end{eqnarray*}
\item induced exterior stacks, i.e.~stacks together with at least one secondary
segment on either or both its sides,
\begin{equation*}
\mathcal{N}'_{\tau}=\mathcal{K}'_{\tau}\times
\left( \mathcal{T}_{\sigma}^2 -1 \right),
\end{equation*}
having generating function
\begin{equation*}
\mathbf{N}'_\tau =\frac{(x y z)^{\tau}}{1-x y z}
\left( {\bf T}_{\sigma}(x){\bf T}_{\sigma}(y)-1\right),
\end{equation*}
\item exterior stems, that is pairs consisting of exterior stacks
$\mathcal{K}'_\tau$ and an arbitrarily long sequence of induced exterior stacks
\begin{equation*}
\mathcal{M}'_\tau=\mathcal{K}'_{\tau} \times
                                \textsc{Seq}\left(\mathcal{N}'_{\tau}\right),
\end{equation*}
having the generating function
\begin{eqnarray*}
\mathbf{M}'_\tau=\frac{\mathbf{K}'_\tau}{1-\mathbf{N}'_\tau}=
\frac{\frac{(x y z)^{\tau}}{1-x y z}}
{1-\frac{(x y z)^{\tau}}{1-x y z}\left( {\bf T}_{\sigma}(x){\bf T}_{\sigma}(y)-1\right)}.
\end{eqnarray*}
\end{itemize}
We inflate all the exterior arcs and the corresponding generating
function is
\begin{eqnarray*}
({\bf M}'_{\tau})^h.
\end{eqnarray*}

{\bf Step III:} here we insert additional secondary segments at the
remaining $(2 t_1+h+1)$ positions in the top and the $(2
t_2+h+1)$ positions in the bottom, see Fig.~\ref{F:RelateJG-3}.
Formally, the third inflation is expressed via the combinatorial class
\begin{equation*}
(\mathcal{T}_{\sigma})^{2 t_1+h+1}\, (\mathcal{T}_{\sigma})^{2 t_2+h+1},
\end{equation*}
where the corresponding generating function is
\begin{eqnarray*}
{\bf T}_{\sigma}(x)^{2 t_1+h+1} \, {\bf T}_{\sigma}(y)^{2 t_2+h+1}.
\end{eqnarray*}
\begin{figure}
\includegraphics[width=1\textwidth]{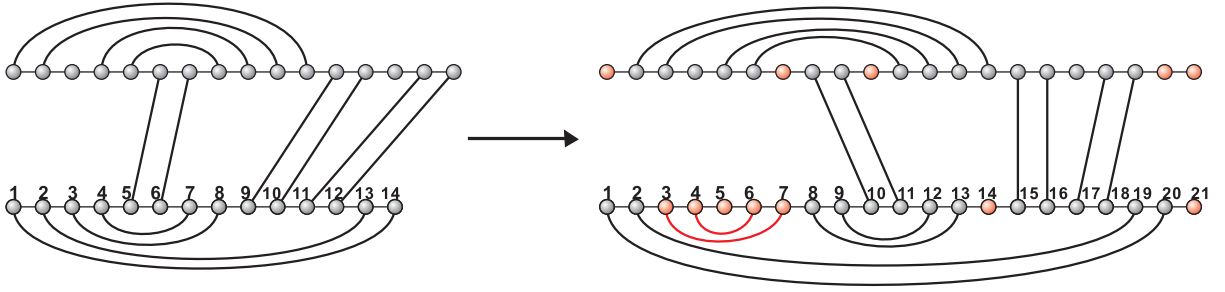}
\caption{Step III: a joint structure (left) obtained in {\sf (1)} in
Fig.~\ref{F:RelateJG-2} is inflated to a new joint structure in ${\mathcal J}_{2,2}$
(right) by adding secondary segments (red).} \label{F:RelateJG-3}
\end{figure}
Combining Step I, Step II and Step III we arrive at
\begin{eqnarray*}
{\mathcal M}_{\sigma}(x)^{t_1}\,{\mathcal M}_{\sigma}(y)^{t_2} \,
({\mathcal M}'_{\tau})^h \,{\mathcal T}_{\sigma}(x)^{2 t_1+h+1} \,
{\mathcal T}_{\sigma}(y)^{2 t_2+h+1}
\end{eqnarray*}
and accordingly
\begin{eqnarray*}
&&{\bf M}_{\sigma}(x)^{t_1}\,{\bf M}_{\sigma}(y)^{t_2} \, ({\bf M}'_{\tau})^h
\,{\bf T}_{\sigma}(x)^{2 t_1+h+1} \, {\bf T}_{\sigma}(y)^{2 t_2+h+1}\\
&=&{\bf T}_{\sigma}(x)\, {\bf T}_{\sigma}(y) ({\bf T}_{\sigma}(x)^2
{\bf M}_{\sigma}(x))^{t_1} ({\bf T}_{\sigma}(y)^2 {\bf M}_{\sigma}(y))^{t_2}
({\bf M}'_{\tau}{\bf T}_{\sigma}(x)\, {\bf T}_{\sigma}(y))^h .
\end{eqnarray*}
Therefore,
\begin{eqnarray*}
 \mathbf{J}_\gamma(x,y,z)  &=&
 {\bf T}_{\sigma}(x)\, {\bf T}_{\sigma}(y) ({\bf T}_{\sigma}(x)^2
 {\bf M}_{\sigma}(x))^{t_1} \\
 &&({\bf T}_{\sigma}(y)^2 {\bf M}_{\sigma}(y))^{t_2}
({\bf M}'_{\tau}{\bf T}_{\sigma}(x)\, {\bf T}_{\sigma}(y))^h.
\end{eqnarray*}
Since for any $\gamma,\gamma_1\in {\mathcal G}(t_1,t_2,h)$ we have
$\mathbf{J}_\gamma(x,y,z)=\mathbf{J}_{\gamma_1}(x,y,z)$, we derive
\begin{equation*}
{\bf J}_{\sigma,\tau}(x,y,z) = \sum_{\gamma\in\,{\mathcal G}}
\mathbf{J}_\gamma(x,y,z) =
\sum_{(t_1,t_2,h) \atop \gamma\in\,{\mathcal G}(t_1,t_2,h)}
G(t_1,t_2,h)\,\mathbf{J}_\gamma(x,y,z).
\end{equation*}
Set
\begin{eqnarray*}
\eta(w)     &=&  \frac{w^{2\sigma}\,{\bf T}_{\sigma}(w)^2}{1-w^2-w^{2\sigma}
({\bf T}_{\sigma}(w)^2-1)}\\
\eta_0      &=&  \frac{(x y z)^{\tau}{\bf T}_{\sigma}(x)\, {\bf T}_{\sigma}(y)}{1- x y z -(x y z)^{\tau}({\bf T}_{\sigma}(x) {\bf T}_{\sigma}(y)-1)}.
\end{eqnarray*}
According to the generating function
\[
{\bf G}(u,v,z)= \sum G(t_1,t_2,h) u^{t_1} v^{t_2} z^h ,
\]
we have
\begin{equation*}
{\bf J}_{\sigma,\tau}(x,y,z) = {\bf T}_{\sigma}(x)\,{\bf T}_{\sigma}(y) \,{\bf G}\!\big(\eta(x),\eta(y),\eta_0\big)
\end{equation*}
and the theorem follows.
\end{proof}


\section{Asymptotic analysis}\label{S:asy}


\subsection{The supercritical paradigm}

\begin{figure}
\includegraphics[width=0.7\textwidth]{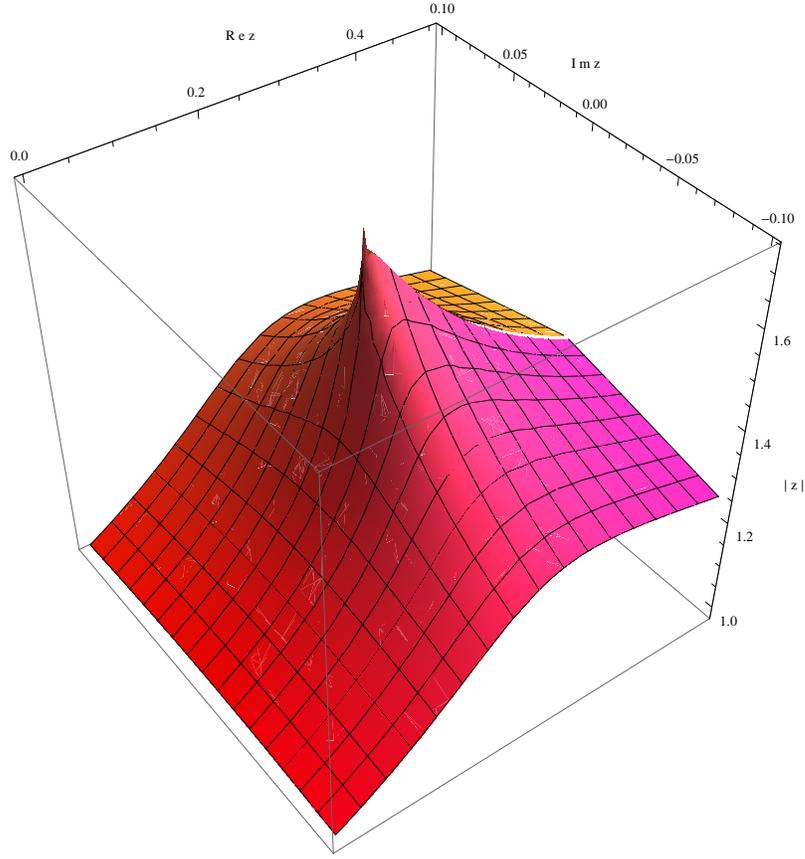}
\caption{Universality of the square root. We display the dominant singularity of
the generating function ${\bf U}(z)$ of shapes (here at $\rho \approx 0.22144$).
All singularities arising from composition of the ``outer'' function
${\bf U}(z)$ governed by the supercritical paradigm produce this type of
singularity, leading to the subexponential factor $n^{-\frac{3}{2}}$.}
\label{F:Singularity}
\end{figure}

Suppose ${\bf U}(z)={\bf G}(z,z,z)$. We view ${\bf U}(z)$ as a generating
function, ${\bf U}(z)= \sum U(l)\,z^l$, where $U(l)$ denotes the number of
shapes having $l$ arcs. It follows from Theorem~\ref{T:Shape} that ${\bf U}
(z)$ satisfies
\begin{equation*}\label{E:U}
(z^2+2z)\,{\bf U}(z)^2-(z^2+3z+1)\, {\bf U}(z) +(1+z)^2=0.
\end{equation*}
Solving this functional equation, we derive
\begin{equation}
{\bf U}(z)= \frac{1+3z+z^2-\sqrt{1-2z-9z^2-10z^3-3z^4}}{2z(z+2)}.
\end{equation}
It is straightforward to verify that the dominant singularity $\rho$ of ${\bf
U}(z)$ is the minimal and positive real solution of
$1-2z-9z^2-10z^3-3z^4=0$ and $\rho \approx 0.22144$, see Fig.~\ref{F:Singularity}.\\
For our computations the following instance of the supercritical
paradigm \citep{Flajolet:07a} is of central importance: we are given a
$D$-finite function, $f(z)$ and an algebraic function $g(u)$
satisfying $g(0)=0$. Furthermore we suppose that $f(g(u))$ has a
unique real valued dominant singularity $\gamma$ and $g$ is regular
in a disc with radius slightly larger than $\gamma$. Then the
supercritical paradigm stipulates that the subexponential factors of
$f(g(u))$ at $u=0$, given that $g(u)$ satisfies certain conditions,
coincide with those of $f(z)$.
\begin{lemma}\label{L:Supercritical}
Let $\vartheta(z)$ be an algebraic, analytic function for $|z|< r$
such that $\vartheta(0)=0$. In addition suppose $\gamma$ is the
unique dominant singularity of ${\bf U}(\vartheta(z))$ and minimum
positive real solution of $\vartheta(z)=\rho$, $|z|< r$,
$\vartheta'(z)\neq 0$. Then ${\bf U}(\vartheta(z))$ has a singular
expansion and
\begin{equation}
[z^n]{\bf U}(\vartheta(z)) \sim c\,n^{-\frac{3}{2}}
\left(\gamma^{-1}\right)^n,
\end{equation}
where $c$ is some constant.
\end{lemma}
\begin{proof}
Since $\vartheta(z)$ is an algebraic function such that
$\vartheta(0)=0$ and ${\bf U}(z)$ is algebraic whence is $D$-finite,
we can conclude that the composition ${\bf U}(\vartheta(z))$ is
$D$-finite. In particular ${\bf U}(\vartheta(z))$ has a singular
expansion.\\
Next, we calculate the singular expansion of the composite function
${\bf U}(\vartheta(z))$. In view of $[z^n]f(z)=\gamma^n
[z^n]f(\frac{z}{\gamma})$ it suffices to analyze the function ${\bf
U}(\vartheta(\gamma z))$ and to subsequently rescale in order to
obtain the correct exponential factor. For this purpose we set
\begin{equation*}
\widetilde{\vartheta}(z)=\vartheta(\gamma z),
\end{equation*}
where $\vartheta(z)$ is analytic in $|z|\leq r$. Consequently
$\widetilde{\vartheta}(z)$ is analytic in $|z|<\widetilde{r}$, for
some $1<\widetilde{r}$. The singular expansion of ${\bf U}(z)$, for
$z\rightarrow \rho$, is given by
\begin{equation*}
{\bf U}(z)                  =   u_0+u_1(\rho -z)^{\frac{1}{2}}( 1 +o(1)).
\end{equation*}
By construction ${\bf U}(\vartheta(\gamma z))={\bf
U}(\widetilde{\vartheta}(z))$, ${\bf U}(\widetilde{\vartheta}(z))$
has the unique dominant singularity at 1. We have the Taylor
expansion of $\widetilde{\vartheta}(z)$ at $z=1$
\begin{equation}\label{E:uinner}
\rho-\widetilde{\vartheta}(z)
=\sum_{n\geq1}\widetilde{\vartheta}_n \,(1-z)^n= \widetilde{\vartheta}_1 \,(1-z)(1+o(1)).
\end{equation}
As for the singular expansion of ${\bf
U}(\widetilde{\vartheta}(z))$, substituting eq.~(\ref{E:uinner})
into the singular expansion of ${\bf U}(z)$, for $z\rightarrow 1$,
\begin{equation*}
{\bf U}(\widetilde{\vartheta}(z))       =   u_0+u_1 \, \widetilde{\vartheta}_1^{\frac{1}{2}} \,(1-z)^{\frac{1}{2}}(1 +o(1)),
\end{equation*}
where $\widetilde{\vartheta}_1=\widetilde{\vartheta}'(z)|_{z=1}
=\gamma \vartheta'(z)|_{z=\gamma} \neq 0$. By
Theorem~\ref{T:transfer1} and Theorem~\ref{T:transfer2} we arrive at
\begin{equation*}
[z^n]{\bf U}(\widetilde{\vartheta}(z))\sim
c\,n^{-\frac{3}{2}}\quad \text{\rm for some constant $c$}.
\end{equation*}
Finally, we use the scaling property of Taylor expansions in order
to derive
\begin{equation*}
[z^n]{\bf U}(\vartheta(z))=\left(\gamma^{-1}\right)^{n}\,[z^n]{\bf
U}(\widetilde{\vartheta}(z))
\end{equation*}
and the proof is complete.
\end{proof}
We remark that Lemma~\ref{L:Supercritical} allows
under certain conditions to obtain the asymptotics of the
coefficients of supercritical compositions of the ``outer'' function
${\bf U}(z)$ and ``inner'' function $\vartheta(z)$. The scenario
considered here is tailored for asymptotic expressions of
$J_{\sigma}(s)$.


\subsection{Asymptotics of $J_{\sigma}(s)$}

\begin{figure*}
\begin{tabular}{cc}
\includegraphics[width=0.5\textwidth]{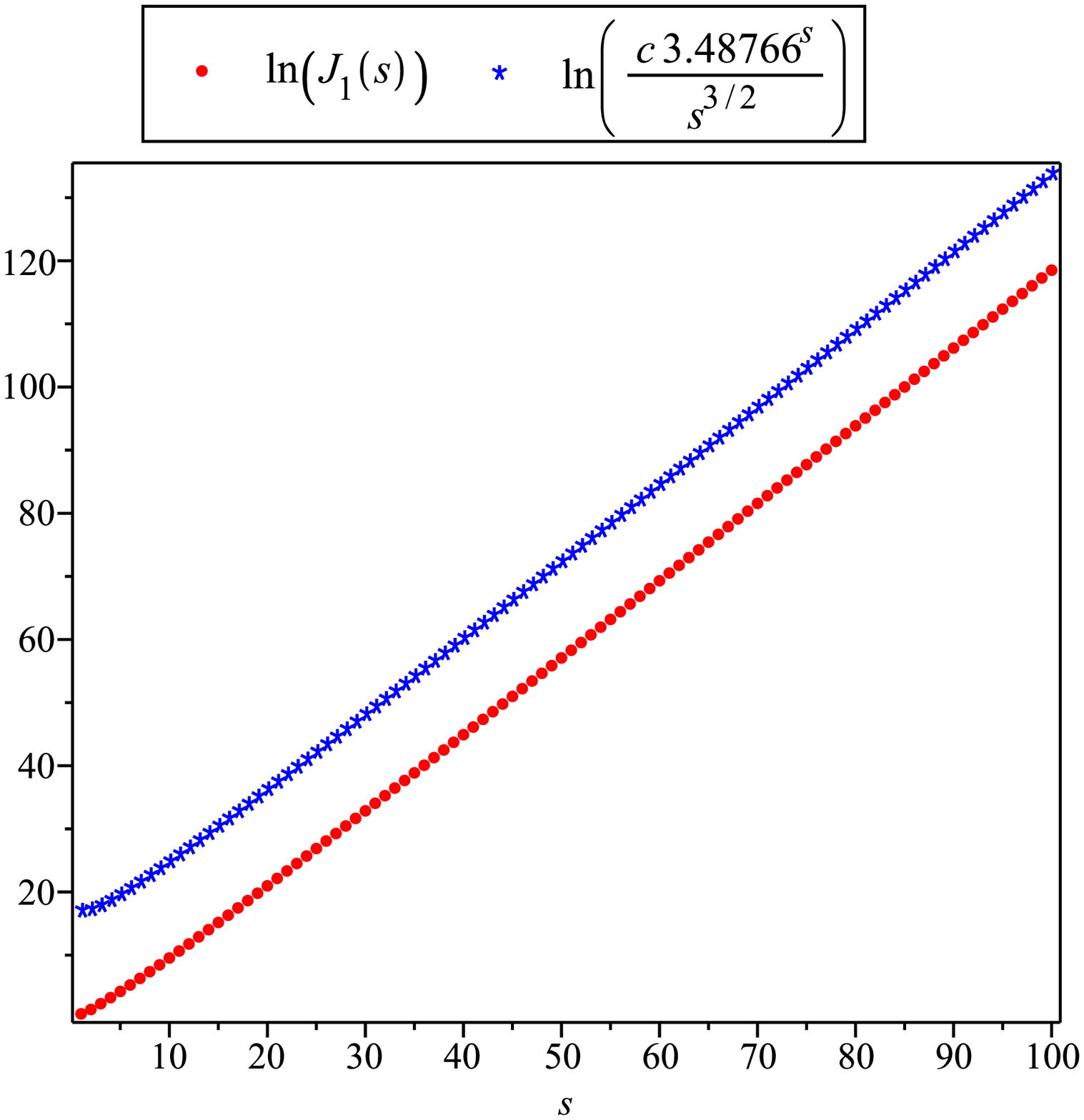}&
\includegraphics[width=0.5\textwidth]{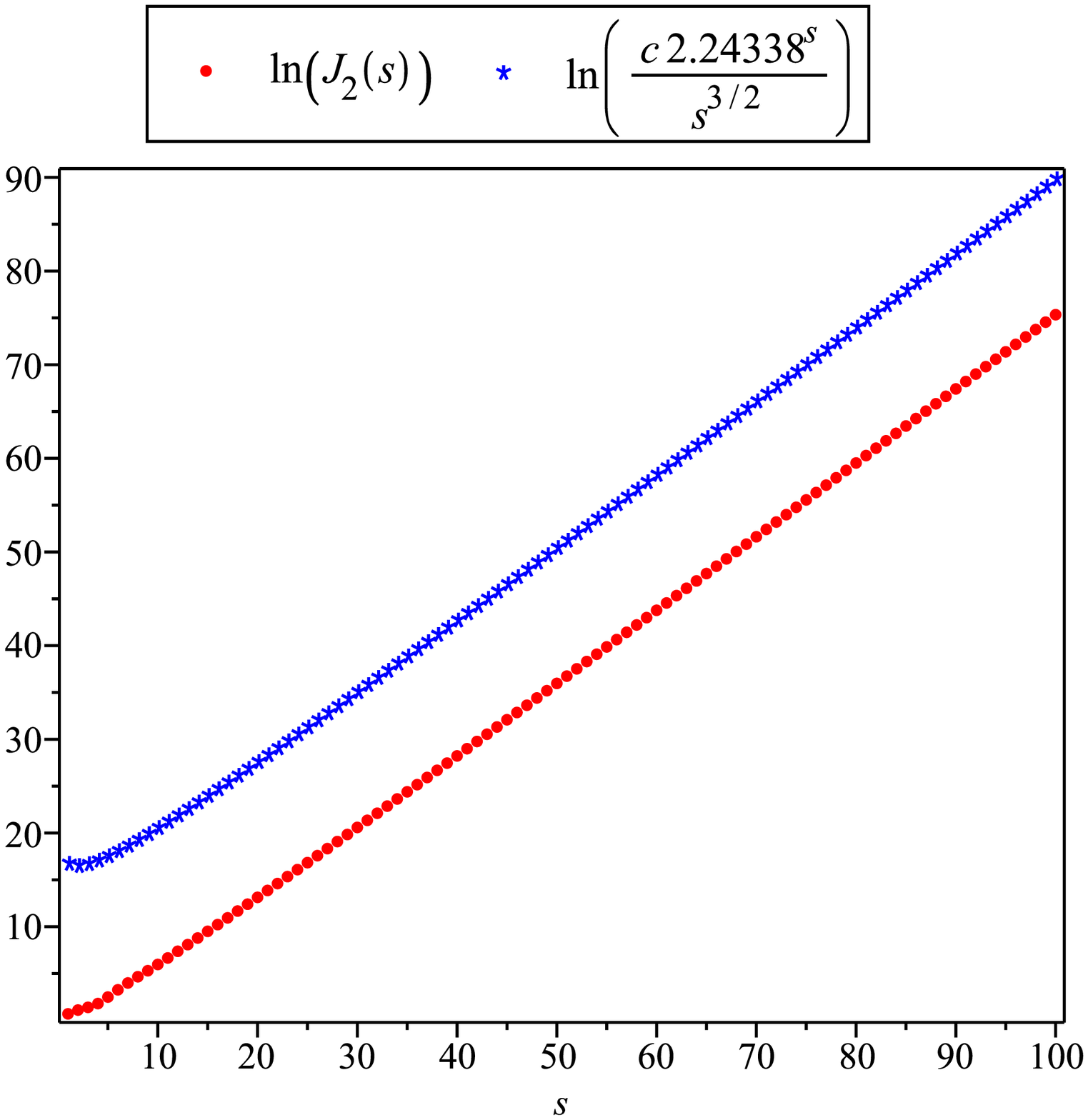}
\end{tabular}
\caption{Exact enumeration versus asymptotic formula. We plot the number
of of joint structures with arc-length $\geq 2$ and stack-length $\geq 1$,
($J_1(s)$) versus its asymptotic formula $c \, s^{-\frac{3}{2}} \,
3.48766^s$ (left) and $J_2(s)$ versus $c \, s^{-\frac{3}{2}} \,
2.24338^s$ (right).
For representational purposes we separate the curves via setting the respective
constants $c=10^7$.
} \label{F:AsymJoint}
\end{figure*}

In this section we shall assume $\sigma=\tau$.
Let $J_{\sigma}(s)$ denote the number of joint
structures of total $s$ vertices having arc-length $\geq 2$,
stack-length $\geq \sigma$ and exterior stack-length $\geq \sigma$
having the generating function
\begin{equation*}
{\bf J}_{\sigma}(z) = \sum J_{\sigma}(s) \, z^s.
\end{equation*}
By definition, we have
\begin{equation*}
{\bf J}_{\sigma}(z) = {\bf J}_{\sigma,\sigma}(z,z,1).
\end{equation*}
\begin{theorem}\label{T:AsymtoticJS}
For $\sigma \geq 1$, we have
\begin{equation}\label{E:JointShape}
{\bf J}_{\sigma}(z) = {\bf T}_{\sigma}(z)^2 \,{\bf U}\!\big(\zeta(z)\big),
\end{equation}
where
\begin{equation}
\zeta(z) = \frac{z^{2\sigma}\,{\bf T}_{\sigma}(z)^2}{1-z^2-z^{2\sigma}(
{\bf T}_{\sigma}(z)^2-1)}.
\end{equation}

Furthermore, for $1 \leq \sigma \leq 9$,
$J_{\sigma}(s)$ satisfies
\begin{equation}\label{E:JointAsymptotic}
J_{\sigma}(s) \, \sim \, c_{\sigma} \, s^{-\frac{3}{2}} \, \left(\gamma_{\sigma}^{-1}\right)^s ,\quad \text{for some
$c_{\sigma}$,}\qquad
\end{equation}
where $\gamma_{\sigma}$ is the minimal, positive real
solution of the equation $\zeta(z)=\rho$, see
Table~\ref{Tab:asymJS}. In particular, $c_{1}\approx
1.6527921$ and $c_{2}\approx 4.3011932$.
\end{theorem}
\begin{proof}
By Theorem~\ref{T:RelateJG} and the definition, we have
\begin{eqnarray*}
{\bf J}_{\sigma}(z) &=& {\bf J}_{\sigma,\sigma}(z,z,1)\\
                    &=& {\bf T}_{\sigma}(z)^2\,{\bf G}\!\big(\zeta(z),\zeta(z),\zeta(z)\big)\\
                    &=& {\bf T}_{\sigma}(z)^2 \,{\bf U}\!\big(\zeta(z)\big),
\end{eqnarray*}
where
\begin{equation*}
\zeta(z) = \frac{z^{2\sigma}\,{\bf T}_{\sigma}(z)^2}{1-z^2-z^{2\sigma}({\bf T}_{\sigma}(z)^2-1)}.
\end{equation*}
Since ${\bf T}_{\sigma}(z)$ is algebraic, we can
conclude that $\zeta(z)$ is algebraic from the closure property of
algebraic functions, whence ${\bf U}(\zeta(z))$ and ${\bf
J}_{\sigma}(z)$ are $D$-finite. Pringsheim's Theorem \citep{Tichmarsh:39}
guarantees that ${\bf J}_{\sigma}(z)$ has a dominant
real positive singularity $\gamma_{\sigma}$. We verify
that for $1\leq \sigma \leq 9$, the
minimal, positive real solution of the equation $\zeta(z)=\rho$ is
strictly smaller than the singularity of $\zeta(z)$, which is
actually the singularity of ${\bf T}_{\sigma}(z)$. Hence
$\gamma_{\sigma}$ is the unique, minimal, positive real
solution of the equation $\zeta(z)=\rho$ and it is straightforward
to check that $\zeta'(z)|_{z=\gamma_{\sigma}}\neq 0$. Therefore the
composite function $ {\bf U}(\zeta(z))$ is governed by the
supercritical paradigm of Lemma~\ref{L:Supercritical}. Furthermore
${\bf T}_{\sigma}(z)$ is analytic at
$\gamma_{\sigma}$, whence the subexponential factors of
${\bf T}_{\sigma}(z)^2 \,{\bf U}(\zeta(z))$ coincide
with those of the function ${\bf U}(z)$. Consequently,
\begin{equation*}
J_{\sigma}(s) \, \sim \, c_{\sigma} \, s^{-\frac{3}{2}} \,
\left(\gamma_{\sigma}^{-1}\right)^s ,\quad \text{for some
$c_{\sigma}$.}\qquad
\end{equation*}
The values of $\gamma_{\sigma}^{-1}$ are listed in Table~\ref{Tab:asymJS}.
It remains to calculate the constant coefficient in the asymptotic formula.
Setting the singular expansion of ${\bf U}(z)$ around $\rho$ and the
Taylor expansions of $\zeta(z)$
and ${\bf T}_{\sigma}(z)^2 $ around $\gamma_{\sigma}$,
\begin{eqnarray*}
{\bf U}(z)    &=&   u_0+u_1(\rho -z)^{\frac{1}{2}}+ O((\rho-z)),\\
\zeta(z)-\rho        &=&   g_1(z-\gamma_{\sigma})+ O((z-\gamma_{\sigma})^2),\\
{\bf T}_{\sigma}(z)^2 &=&   t_0+t_1(\gamma_{\sigma}-z)+O((\gamma_{\sigma}-z)^2).
\end{eqnarray*}
We proceed by substituting these expansions into ${\bf
T}_{\sigma}(z)^2 \,{\bf U}\!\left(\zeta(z)\right)$
\[
{\bf J}_{\sigma}(z)= t_0 u_0 + t_0 u_1 g_1^{\frac{1}{2}}
(\gamma_{\sigma}-z)^{\frac{1}{2}}+ O(\gamma_{\sigma}-z).
\]
Using Theorem~\ref{T:transfer1} and Theorem~\ref{T:transfer2}, we
have
\begin{equation*}
J_{\sigma}(s) \, \sim \, \frac{t_0 u_1
(g_1 \gamma_{\sigma})^{\frac{1}{2}}}{\Gamma (-\frac{1}{2})} \,
s^{-\frac{3}{2}} \, \left(\gamma_{\sigma}\right)^{-s}
\end{equation*}
Setting $c_{\sigma} = \frac{t_0 u_1 (g_1
\gamma_{\sigma})^{\frac{1}{2}}}{\Gamma (-\frac{1}{2})}$, we compute
$c_{1}\approx 1.6527921$ and
$c_{2}\approx 4.3011932$, completing the proof of Theorem~\ref{T:AsymtoticJS}.
\end{proof}
\begin{table}
\begin{center}
\setlength{\tabcolsep}{3.5pt}
\caption{Exponential growth rates
$\gamma_{\sigma}^{-1}$ for joint structures
with arc-length $\geq 2$, having both stack-length and
exterior stack-length $\geq \sigma$.}
\label{Tab:asymJS}
\begin{tabular}{cccccccccc}
\hline\noalign{\smallskip}
\small$\sigma$                &\small$1$                     & \small$2$       & \small$3$        & \small$4$       & \small $5$       & \small$6$       &\small $7$        & \small $8$       & \small $9$  \\
\noalign{\smallskip}\hline\noalign{\smallskip}
\small$\gamma_{\sigma}^{-1}$  &\small$3.48766$               & \small$2.24338$ & \small $1.86724$ &\small $1.67974$ & \small $1.56544$ &\small $1.48763$ & \small $1.43083$ & \small $1.38731$ & \small $1.35276$ \\
\noalign{\smallskip}\hline
\end{tabular}
\end{center}
\end{table}
We next observe that eq.~(\ref{E:JointShape}) allows us to derive a functional
equation for ${\bf J}_{\sigma}(z)$, which in turn gives a
recurrence of $J_{\sigma}(s)$.
\begin{corollary}\label{C:EnumerateJoint}
For $\sigma \geq 1$, the generating function ${\bf J}_{\sigma}(z)$ satisfies
the functional equation
\begin{equation}\label{E:JointEq}
{\bf A}(z){\bf J}_{\sigma}(z)^2+ {\bf B}(z){\bf J}_{\sigma}(z)+{\bf C}(z)=0,
\end{equation}
where
\begin{equation}
\begin{aligned}
{\bf A}(z)  = &z^{2\sigma}\,(2-2z+2z^{2\sigma}-z^{2\sigma} {\bf T}_{\sigma}(z)^2 ),\\
{\bf B}(z)  = &-\Big(1-2z^2+z^4+(2+{\bf T}_{\sigma}(z)^2) z^{2\sigma}\\
              & -(2+{\bf T}_{\sigma}(z)^2) z^{2+2\sigma}+
(1+{\bf T}_{\sigma}(z)^2-{\bf T}_{\sigma}(z)^4) z^{4\sigma}\Big),\\
{\bf C}(z)  = &(1-z^2+z^{2\sigma})^2.
\end{aligned}
\end{equation}
Furthermore, the number $J_{\sigma}(s)$ of joint structures with total $s$ vertices
satisfies the following recurrence:
\begin{equation*}
J_{\sigma}(s)= c(s) + \sum_{i=1}^{s} b(i)\, J_{\sigma}(s-i) +\sum_{i=1}^{s}
\sum_{j=0}^{s-i} a(i) \, J_{\sigma}(j)\,J_{\sigma}(s-i-j),
\end{equation*}
where $a(s)$, $b(s)$ and $c(s)$ are the coefficients of $z^s$ of ${\bf A}(z)$,
${\bf B}(z)$ and ${\bf C}(z)$, respectively.
\end{corollary}
In Table~\ref{Tab:JointEnumerate}, we list the numbers of joint
structures $J_1(s)$ and $J_2(s)$ for $s=1,\dots, 12$.
\begin{proof}
Substituting
$z=\frac{z^{2\sigma}\,{\bf T}_{\sigma}(z)^2}{1-z^2-z^{2\sigma}
({\bf T}_{\sigma}(z)^2-1)}$ into eq.~(\ref{E:U}) and using
eq.~(\ref{E:JointShape}), we obtain eq.~(\ref{E:JointEq}).
Note that $a(0)=0$ and $b(0)=-1$. Calculating the
coefficients of $z^s$ of eq.~(\ref{E:JointEq}), the recurrence
follows immediately.
\end{proof}
\begin{table}
\caption{\small The numbers of joint structures $J_1(s)$ and $J_2(s)$ over a total number
of $s=1,\dots, 12$ nucleotides.} \label{Tab:JointEnumerate}
\begin{center}
\begin{tabular}{c*{12}{c}}
\hline\noalign{\smallskip}
$s$     & \small$1$ & \small$2$ & \small$3$ & \small$4$ & \small $5$ &
\small$6$ &\small $7$ & \small $8$ & \small $9$ & \small $10$ & \small
$11$ & \small $12$ \\
\noalign{\smallskip}\hline\noalign{\smallskip}
$J_{1}(s)$     & \small$2$ & \small$4$ & \small$10$ & \small$26$ &
\small $70$ & \small$194$ &\small $550$ & \small $1590$ & \small
$4674$ & \small $13940$ & \small $42106$ & \small $128610$\\
$J_{2}(s)$     & \small$2$ & \small$3$ & \small$4$  & \small$6$ &
\small $12$ & \small$26$ &\small $54$ & \small $105$ & \small
$200$ & \small $389$ & \small $780$ & \small $1589$\\
\noalign{\smallskip}\hline
\end{tabular}
\end{center}
\end{table}
In Fig.~\ref{F:AsymJoint}, we show that our asymptotic formulas
work well already for small sequence length.
Here we contrast the exact values, $J_1(s)$ and $J_2(s)$, with the
asymptotic formulas given via Theorem~\ref{T:AsymtoticJS}:
\begin{equation*}
J_1(s) \, \sim \, c_1 \, s^{-\frac{3}{2}} \, 3.48766^s \quad \text{and}
\quad J_2(s) \, \sim \, c_2 \, s^{-\frac{3}{2}} \, 2.24338^s .
\end{equation*}

\section{Discussion}\label{S:discussion}


The discovery of more and more instances of regulatory actions among
RNA molecules make evident that RNA-RNA interaction is a problem of
central importance.
While it is wellknown how to MFE-fold these interaction structures
\citep{Alkan:06,rip:09} this paper constitutes progress with respect
to the theoretical understanding of RNA-RNA interaction structures.
Insights in the combinatorics of joint structures allows deeper
understanding of analysis and design of folding algorithms as well as
algorithmic approximations.

At first sight, it should be straightforward to derive the generating
function of joint structures from the (eleven) recursion relations of
the original {\tt rip}-grammar (implied by Proposition~\ref{P:Decomposition})
\citep{rip:09}.
While this is {\it in principle} correct, the mere statement of the
generating function derived this way fills several pages. This approach is
neither suitable
for deriving any asymptotic formulas nor does it allow to deal with specific
stack-length conditions. In fact, the extraction of its coefficients would present a
nontrivial task.

We do not use the recurrences of \citep{rip:09} directly.
Instead we build our theory of
joint structures centered around the concept of shapes.
The key to all results is the simple shape-grammar of Theorem~\ref{T:Shape}.
The basic idea here is that the collapsing of stems preserves vital information
of the interaction structure. Given a shape a joint structure can be
obtained via inflation, see Theorem~\ref{T:RelateJG}.

While there exists a notion of shapes for RNA secondary structures \citep{Giegerich}
their combinatorics is not shape-based. Everything is organized around
recurrences, which oftentimes hides deeper structural insight and connections.
As a result symbolic enumeration has not been employed in order to derive the
generating function of RNA secondary structures.

In contrast, RNA pseudoknot structures  \citep{Modular} represent a
shape-based structure class (here further complication enters the picture as the
generating function of their shapes can only be computed via the reflection
principle).

The theory of joint structures presented here resembles features of
the theory of modular diagrams and is in particular shape-based. However, the shapes
of joint structures are governed by simple algebraic generating functions and satisfy
a simple recurrence.

Let us finally outline future research. We currently study the generating
function of canonical joint structures having minimum arc-length four. This
derivation requires a more detailed look at shapes of joint structures since
additional variables have to be introduced. The purpose of these variables is to
allow to distinguish specific inflation scenarious.


\section{Further results}\label{S:later}

In this section we generalize our results to joint structures with arc-length
$\geq \lambda$, interior stack-length $\geq \sigma$, and exterior stack-length
$\geq \tau$.
Let $\mathcal{T}_{\sigma}^{[\lambda]}$ denote the combinatorial
class of $\sigma$-canonical secondary structures having arc-length
$\geq \lambda$ and $T_{\sigma}^{[\lambda]}(n)$ denote the number of
all $\sigma$-canonical secondary structures with $n$ vertices having
arc-length $\geq \lambda$ and
\begin{equation*}
{\bf T}^{[\lambda]}_{\sigma}(z)= \sum T_{\sigma}^{[\lambda]}(n)\, z^n.
\end{equation*}
\begin{theorem}\label{T:SecondSigLamb}
Let $\sigma\in\mathbb{N}$, $z$ be an indeterminant and let
\begin{eqnarray*}
u_\sigma(z) & = & {(z^2)^{\sigma-1}\over{z^{2\sigma}-z^2+1}}, \\
v_\lambda(z) & = & 1-z+u_\sigma(z) {\sum_{h=2}^\lambda z^h},
\end{eqnarray*}
then, ${\bf T}^{[\lambda]}_{\sigma}(z)$, the generating function of
$\sigma$-canonical structures with minimum arc-length $\lambda$ is given by
\begin{equation*}
{\bf T}_{\sigma}^{[\lambda]}(z) = \frac{1}{v_\lambda(z)} {\bf
F}\left(\left(\frac{\sqrt{u_\sigma(z)}\,z}{v_\lambda(z)}
\right)^2\right),
\end{equation*}
where
\begin{equation*}
{\bf F}(z)   = \frac{1-\sqrt{1-4z}}{2z}.
\end{equation*}
\end{theorem}
Theorem~\ref{T:SecondSigLamb} implies that
${\bf T}^{[\lambda]}_{\sigma}(z)$ is an algebraic function for any
specified $\lambda$ and $\sigma$, since ${\bf F}(z)$ is algebraic
and $v_\lambda(z), u_\sigma(z)$ are both rational functions.

We are now in position to establish a generalization
of Theorem~\ref{T:RelateJG} that allows us to compute the generating
function ${\bf J}_{\sigma,\tau}^{[\lambda]}(x,y,z)$ for $ \lambda
\leq \tau+1$.\\
Let ${\mathcal J}_{\sigma,\tau}^{[\lambda]}$ denote the class of joint
structures with arc-length $\geq \lambda$, interior stack-length $\geq
\sigma$, and exterior stack-length $\geq \tau$. Let
$J_{\sigma,\tau}^{[\lambda]}(n,m,h)$ denote the number of joint structures
in ${\mathcal J}_{\sigma,\tau}^{[\lambda]}$ having $n$ vertices in the top,
$m$ vertices in the bottom, $h$ exterior arcs having the generating
function
\begin{equation*}
{\bf J}_{\sigma,\tau}^{[\lambda]}(x,y,z) =
\sum J_{\sigma,\tau}^{[\lambda]}(n,m,h)\, x^{n} y^{m} z^h.
\end{equation*}
\begin{theorem}\label{T:RelateJGLamb}
For $\sigma \geq 1,\tau \geq 1,\, \lambda \leq \tau+1$ , we have
\begin{equation}
{\bf J}_{\sigma,\tau}^{[\lambda]}(x,y,z) =
{\bf T}_{\sigma}^{[\lambda]}(x)\,{\bf T}_{\sigma}^{[\lambda]}(y) \,
{\bf G}\!\big(\eta(x),\eta(y),\eta_0\big),
\end{equation}
where
\begin{eqnarray*}
\eta(w)     &=&  \frac{w^{2\sigma}\,
{\bf T}_{\sigma}^{[\lambda]}(w)^2}{1-w^2-w^{2\sigma}
({\bf T}_{\sigma}^{[\lambda]}(w)^2-1)}\\
\eta_0      &=&  \frac{(x y z)^{\tau}{\bf T}_{\sigma}^{[\lambda]}(x)\,
{\bf T}_{\sigma}^{[\lambda]}(y)}{1- x y z -(x y z)^{\tau}({\bf T}_{\sigma}^{[
\lambda]}(x) {\bf T}_{\sigma}^{[\lambda]}(y)-1)}.
\end{eqnarray*}
\end{theorem}
\begin{proof}
Using the notation and approach of Theorem~\ref{T:RelateJG} we
arrives at
\begin{align*}
\mathcal{K}_{\sigma}    &=  \mathcal{R}^{\sigma}\times\textsc{Seq}\left(\mathcal{R}\right)\\
\mathcal{N}_{\sigma}    &=  \mathcal{K}_{\sigma}\times\left( (\mathcal{T}_{\sigma}^{[\lambda]})^2 -1 \right)\\
\mathcal{M}_{\sigma}    &=  \mathcal{K}_{\sigma} \times \textsc{Seq}\left(\mathcal{N}_{\sigma}\right)\\
\mathcal{K}'_{\tau}    &=  \mathcal{R}_0^{\tau}\times\textsc{Seq}\left(\mathcal{R}_0\right)\\
\mathcal{N}'_{\tau}    &=  \mathcal{K}'_{\tau}\times\left( (\mathcal{T}_{\sigma}^{[\lambda]})^2 -1 \right)\\
\mathcal{M}'_{\tau}    &=  \mathcal{K}'_{\tau} \times \textsc{Seq}\left(\mathcal{N}'_{\tau}\right)\\
\mathcal{J}_{\sigma,\tau}^{[\lambda]}   &=  {\mathcal M}_{\sigma}(x)^{t_1}\,{\mathcal M}_{\sigma}(y)^{t_2}\,
({\mathcal M}'_{\tau})^h \,(\mathcal{T}_{\sigma}^{[\lambda]}(x))^{2 t_1+h+1} \, (\mathcal{T}_{\sigma}^{[\lambda]}(y))^{2 t_2+h+1}.
\end{align*}
The only difference is that $\mathcal{T}_{\sigma}^{[\lambda]}$
replaces $\mathcal{T}_{\sigma}$ to make the structure with
arc-length $\geq \lambda$. The key point here is that the
restriction $\lambda \leq \tau+1$ guarantees that any $2$-arc in
$\gamma$ has after inflation a minimum arc-length of
$\tau+1\geq\lambda$.\\
Therefore, the generating function of class
$\mathcal{J}_{\sigma,\tau}^{[\lambda]}$ satisfies
\begin{equation*}
{\bf J}_{\sigma,\tau}^{[\lambda]}(x,y,z) = {\bf T}_{\sigma}^{[\lambda]}(x)\,{\bf T}_{\sigma}^{[\lambda]}(y) \,{\bf G}\!\big(\eta(x),\eta(y),\eta_0\big),
\end{equation*}
where
\begin{eqnarray*}
\eta(w)     &=&  \frac{w^{2\sigma}\,{\bf T}_{\sigma}^{[\lambda]}(w)^2}{1-w^2-w^{2\sigma} ({\bf T}_{\sigma}^{[\lambda]}(w)^2-1)}\\
\eta_0      &=&  \frac{(x y z)^{\tau}{\bf T}_{\sigma}^{[\lambda]}(x)\, {\bf T}_{\sigma}^{[\lambda]}(y)}
{1- x y z -(x y z)^{\tau}({\bf T}_{\sigma}^{[\lambda]}(x) {\bf T}_{\sigma}^{[\lambda]}(y)-1)}.
\end{eqnarray*}
\end{proof}
We remark that Theorem~\ref{T:RelateJGLamb} immediately implies
Theorem~\ref{T:RelateJG}.\\
Analogously, we have
\begin{theorem}
For $\lambda \leq \sigma +1$, we have
\begin{equation}
{\bf J}_{\sigma}^{[\lambda]}(z) = {\bf T}_{\sigma}^{[\lambda]}(z)^2 \,{\bf U}\!\big(\zeta(z)\big),
\end{equation}
where
\begin{equation}
\zeta(z) = \frac{z^{2\sigma}\,{\bf T}_{\sigma}^{[\lambda]}(z)^2}{1-z^2-z^{2\sigma}({\bf T}_{\sigma}^{[\lambda]}(z)^2-1)}.
\end{equation}

Furthermore, for $1 \leq \sigma \leq 9$ and $1 \leq \lambda \leq 5$,
$J_{\sigma}^{[\lambda]}(s)$ satisfies
\begin{equation}\label{E:JointAsymptoticLamb}
J_{\sigma}^{[\lambda]}(s) \, \sim \, c_{\sigma}^{[\lambda]} \, s^{-\frac{3}{2}} \, \left(\frac{1}{\gamma_{\sigma}^{[\lambda]}}\right)^s ,\quad \text{for some
$c_{\sigma}^{[\lambda]}$,}\qquad
\end{equation}
where $\gamma_{\sigma}^{[\lambda]}$ is the minimal, positive real
solution of the equation $\zeta(z)=\rho$, see
Table~\ref{Tab:asymJSLamb}. In particular, $c_{1}^{[2]}\approx
1.6527921$, $c_{2}^{[2]}\approx 4.3011932$, and $c_{2}^{[3]}\approx
3.8671841$.
\end{theorem}

\begin{table}
\setlength{\tabcolsep}{3.5pt}
\caption{Exponential growth rates
$\left(\gamma_{\sigma}^{[\lambda]}\right)^{-1}$ for joint structures
with arc-length $\geq \lambda$, having both stack-length and
exterior stack-length $\geq \sigma$.}
\label{Tab:asymJSLamb}
\begin{center}
\begin{tabular}{cccccccccc}
\hline\noalign{\smallskip}
\small$\sigma$     &\small$1$                     & \small$2$       & \small$3$        & \small$4$       & \small $5$       & \small$6$       &\small $7$        & \small $8$       & \small $9$  \\
\noalign{\smallskip}\hline\noalign{\smallskip}
\small$\lambda=1$  &\small$3.77438$               & \small$2.30663$ & \small $1.89559$ &\small $1.69615$ & \small $1.57629$ &\small $1.49541$ & \small $1.43671$ & \small $1.39194$ & \small $1.35651$ \\
\small$\lambda=2$  &\small$3.48766$               & \small$2.24338$ & \small $1.86724$ &\small $1.67974$ & \small $1.56544$ &\small $1.48763$ & \small $1.43083$ & \small $1.38731$ & \small $1.35276$ \\
\small$\lambda=3$  &\small$0.00000$               & \small$2.21090$ & \small $1.84998$ &\small $1.66876$ & \small $1.55773$ &\small $1.48187$ & \small $1.42633$ & \small $1.38368$ & \small $1.34976$ \\
\small$\lambda=4$  &\small$0.00000$               & \small$0.00000$ & \small $1.83971$ &\small $1.66155$ & \small $1.55233$ &\small $1.47764$ & \small $1.42291$ & \small $1.38085$ & \small $1.34737$ \\
\small$\lambda=5$  &\small$0.00000$               & \small$0.00000$ & \small $0.00000$ &\small $1.65691$ & \small $1.54861$ &\small $1.47459$ & \small $1.42036$ & \small $1.37867$ & \small $1.34549$ \\
\noalign{\smallskip}\hline
\end{tabular}
\end{center}
\end{table}
{\bf Acknowledgments.}
We would like to thank F.W.D.\ Huang, E.Y. Jin and R.R.\ Wang for
discussions.

%
%

\bibliographystyle{spbasic}

\end{document}